\documentclass[12pt]{article}
\usepackage{amsmath}
\usepackage{amsfonts}

\newtheorem{remark}{Remark}
\newtheorem{theorem}{Theorem}

\newtheorem{lemma}{Lemma}
\newtheorem{definition}{Definition}[section]
\newtheorem{corollary}{Corollary}
\newcommand{\riangle}{\right\rangle}
\newcommand{\eq}{\begin{equation}\begin{array}{lllllllllllllllllllllllllllllllll}}
\newcommand{\ee}{\end{array}\end{equation}}
\newcommand{\bmt}{\left( \begin{array}{ccccccccc}}
\newcommand{\emt}{\end{array}\right)}
\newcommand{\bea}{\begin{eqnarray}}
\newcommand{\eea}{\end{eqnarray}}
\newcommand{\bean}{\begin{eqnarray*}}
\newcommand{\eean}{\end{eqnarray*}}
\newcommand{\leangle}{\left\langle}

\newcommand{\U}{{\bf U}}

\newcommand{\prf}{{\bf Proof.\quad }}
\newcommand{\qed}{\hfill\hbox{\hskip 1pt \vrule width 7pt height 6pt 
			depth 1.5pt \hskip 1pt}}
\newenvironment{pmx}[1]{\left(\begin{array}{#1}}{\end{array}\right)}
\newcommand{\cD}{\mathcal{D}}

\newcommand{\HI}{\mathbb{H}}

\newcommand{\fb}{\mathfrak{b}}
\newcommand{\fg}{\mathfrak{g}}

\newcommand{\fh}{\mathfrak{h}}
\newcommand{\fp}{\mathfrak{p}}
\newcommand{\fk}{\mathfrak{k}}

\newcommand{\half}{{1\over2}}

\newcommand{\Det}{{\rm Det\,}}

\newcommand{\vp}{{\varphi}}

\newcommand{\IR}{{\mathbb R}}

\begin{document}

\title{The Symplectic Structure of Curves in three dimensional spaces 
of constant curvature and the equations of mathematical physics}
\author{V. Jurdjevic} 
\date{University of Toronto}
\maketitle
\begin{center}

\end{center}
\section{Introduction}

    This paper  defines a symplectic form on the infinite dimensional Fr\'echet manifold of framed curves of fixed length over a three dimensional  simply connected Riemannian manifold of constant curvature. The framed curves are anchored at the initial point and are further constrained by the condition that the tangent vector of the projected curve coincides with the first leg of the orthonormal frame. Such class of curves are called  anchored Darboux curves and in particular include the Serret-Frenet framed curves.

The symplectic form $\omega$ is defined on the universal covers of the orthonormal frame bundles of the underlying manifolds: $SL_2(C)$ for the hyperboloid $\HI^3$, $SU_2\times SU_2$ for the sphere $S^3$, and the semidirect product $E^3\triangleright SU_2$ for the Euclidean space $E^3$. The form $\omega$ is left invariant on each of the above groups, and is  induced by the Poisson-Lie bracket on the  appropriate Lie algebra. More precisely, the form $\omega$ in each of the above non-Euclidean cases  is defined on the Cartan space $\fp$  corresponding to the decomposition \[\fg=\fp+\fk \] of the Lie algebra $\fg$   subject  to the usual Lie algebraic relations\[[\fp,\fp]=\fk,\,[\fp,\fk]=\fp,\,[\fk.\fk]=\fk\]
   In the case of the hyperboloid $\fg$ is equal to $sl_2(C)$ and the Cartan space is equal to the  space of  the Hermitian matrices $\fh$, and in the case of the sphere $\fg$ is equal to $su_2\times su_2$ and the  Cartan space is isomorphic to the space of skew-Hermitian matrices $\fk$.  In each case the maximal compact subalgebra $\fk$ is equal to the Lie algebra of skew Hermitian matrices. The symplectic forms in each of these two cases are isomorphic  to each other as a consequence of the isomorphism  between $\fk$ and $\fh$ given by $i\fh=\fk$. The Euclidean space $E^3$ is identified with $\fh$ equipped with the metric defined by the trace form, and its framed curves are represented in the semidirect product $\fh\triangleright \fk$. The Euclidean Darboux curves inherit the hyperbolic symplectic form $\omega$ which is isomorphic to the symplectic form used by J. Millson and B. Zombro in ([17]).

 Each group $G$ mentioned above is a principal $SU_2$ bundle over the underlying symmetric space with a natural connection defined by the left invariant vector fields that  take values in the Cartan space $\fp$. The vertical distribution is defined by the left invariant vector fields that take their values in $\fk$.  
  In this setting then, anchored Darboux curves are the solutions $g(s)\in G$ of a differential equation\eq\frac{dg}{ds}(s)=g(s)(E_1+u_1(s)A_1+u(s)A_2+u_3(s)A_3)\ee
with $g(0)=I$, where $E_1$ is a fixed unit vector in the Cartan space $\fp$. The matrices $A_1,A_2,A_3$ denote the  skew-Hermitian Pauli matrices, and $u_1(s),u_2(s),u_3(s)$ are arbitrary  real valued functions
on a fixed interval $[0,L]$.  Each anchored Darboux curve defines a horizontal Darboux curve $h(s)\in G$ that is a solution of the differential equation \eq\label{lambda}\frac{dh}{ds}(s)=h(s)\Lambda (s)\mbox{, }\Lambda (s)=R(s)E_1R^{-1}(s)\ee  with $R(s)$ the solution curve in $SU_2$ of  the equation \eq\label{R}\frac{dR}{ds}=R(s)(u_1(s)A_1+u_2(s)A_2+u_3(s)A_3)\ee that satisfies $R(0)=I$.
 The symplectic form for the hyperbolic Darboux curves is given by\eq\omega_\Lambda (V_1,V_2)={1\over i}\int_0^L\langle \Lambda (s),[U_1(s),U_2(s)]\rangle\,ds\ee with $ U_1(s)$ and $U_2(s)$ Hermitian matrices orthogonal to the tangent vector $\Lambda (s)$,
 that further satisfy $U_j(0)=0$ and $\frac{dV_j}{ds}(s)=U_j(s)$ for $j=1,2$.

In the spherical case the symplectic form  has the same form as in the hyperbolic case, except for the factor $1\over i$, which is omitted. The matrices $U_j$  in this case take values in $\fk$ and satisfy \[\frac{dV_j}{ds}(s)=[\Lambda (s),V_j(s)]+U_j(s)\] for $j=1,2$.\\

The second part of the paper is devoted to the Hamiltonian flow associated with the function
\[f(g(s))=\half\int_0^L||\frac{d\Lambda }{ds}(s)||^2\,ds=\half\int_0^L\kappa^2(s)\,ds\] where $g$  denotes a frame-periodic horizontal Darboux periodic  curve, i.e., a Darboux curve for which the solution $R(s)$ of equation (\ref{R}) is periodic . Here $\kappa (s)$ denotes the curvature of the projected curve $x(s)$ in the underlying  symmetric space. 
 
On the level of Lie algebras, the Hamiltonian flow induced by the symplectic form $\omega$   generates Heisenberg's magnetic equation in the Cartan space $\fp$ given by \eq
  \frac{\partial \Lambda}{\partial t}(s,t)={1\over i}[\Lambda (s),\frac{\partial ^2\Lambda}{\partial s^2}(s,t)]\ee in the hyperbolic and the Euclidean case, and by\[\frac{\partial \Lambda}{\partial t}(s,t)=[\Lambda (s),\frac{\partial ^2\Lambda}{\partial s^2}(s,t)]\] in the spherical case. 

The corresponding frame $R(s,t)$ defines a complex function \eq\label{psi}
\psi (s,t)=u(s,t)\exp({i\int_0^su_1(x,t)\,dx})\ee with $u(s,t)=u_2(s,t)+iu_3(s,t)$ that is a solution of the non-linear Schroedinger's equation \eq\label{sch}
-i\frac{\partial \psi}{\partial t}(t,s)=\frac{\partial^2 \psi}{\partial s^2}(t,s)+1/2|\psi(t,s)|^2(\psi(t,s)+c)
\ee
(Theorem \ref{h-s}).

This finding clarifies a remarkable observation of H. Hasimoto ([8]) that  the function\[\psi (s,t)=\kappa (s,t)\exp({i\int_0^s\tau (x,t)\,dx})\]
where $\kappa(t,s)$ and $\tau(t,s)$ are the curvature and the torsion of a
curve $\gamma(t,s))$ that evolves according to to the filament equation
\begin{equation}
\frac{\partial \gamma}{\partial t}(t,s)=\kappa(t,s)B(t,s)
\end{equation}
is a solution of the non-linear Schroedinger equation
\[
-i\frac{\partial \psi}{\partial t}(t,s)=\frac{\partial^2 \psi}{\partial s^2}(t,s)+1/2|\psi(t,s)|^2(\psi(t,s)+c)
\]
for some complex constant c. Indeed, when the frame $R(s)$ in equation (\ref{R}) is a Serret-Frenet frame then $\psi$ given by (\ref{psi}) coincides with  Hasimoto's function up to a unitary constant.

The curves that correspond to the critical points of $f=\half\int_0^L\kappa^2(s)\,ds$ are called elastic. The material in Section $5$ shows that the elastic curves  with periodic curvatures always generate  periodic soliton solutions for the non-linear Schroedinger's equation.  The extremal equations associated with elastic curves, obtained through the Hamiltonian formalism of Lie groups, provide for easy transitions to solitons:

The Hamiltonian system associated with elastic curves is completely integrable having four integrals of motion: the Hamiltonian itself, which may be regarded as the energy of the elastic curve, two Casimirs, and another integral of motion  $H_1$ due to  an extra symmetry in the system. The elastic curves that generate solitons reside on a fixed energy level and propagate with the speed equal to $H_1$. The fact that the equations for the heavy top form an invariant subsystem of the equations for the elastic curves ([10] and [11]) makes the connection between elastic curves and solitons  even more intriguing: the  speed of the soliton corresponds to the angular momentum along the axis of symmetry for the top of Lagrange.

The formalism of this paper suggests that there is a class of functions $f_0,f_1,f_2\ldots$   over the space of Darboux curves that begins with $f_0=\half\int_0^L\kappa^2(s)\,ds$ having the property that any two functions Poisson commute. It is shown  in the paper that $f_1$ and  $f_2$ given by
\[f_1=i\int_0^L\langle [\Lambda (s),\frac{d\Lambda}{ds}(s)],\frac{d^2\Lambda}{ds^2}(s)\rangle\,ds,\,f_2=\int_0^L(||\ddot{\Lambda}(s)||^2-\frac{5}{4}||\dot{\Lambda}(s)||^4 )\,ds\]
  are in this class. 

The above funtions can be expressed  either in terms of the geometric invariants of the  underlying Darboux curve as: \[f_1=\int_0^L \kappa^2(s)\tau (s)\,ds,\,f_2=\int_0^L (\frac{\partial \kappa}{\partial s}(s)^2+\kappa^2(s)\tau^2(s)-\frac{1}{4}\kappa^4(s) )\,ds\] in which  case they agree with the first three functions on the list presented by J. Langer and R. Perline in ([14]), or they can be expressed in terms of  the complex function $u(s)$ defined by equation (\ref{psi}) as $f_0=\half\int_0^L|u(s)|^2\,ds$ and \[f_1={1\over 2i}\int_0^L(\bar{u}\dot{u}-{u}\dot{\bar{u}})\,ds,\,f_2=\int_0^L(|\frac{\partial u}{\partial s}(s,t)|^2 -\frac{1}{4}|u(s,t)|^4)\,ds\]
  in which case they correspond to the  first three conserved quantities, the number of particles, the momentum and the energy, in the paper by C.Shabat and V. Zakharov in ([18]).

The paper is organized as follows. The first part of the paper leads up to the symplectic form.
 The fact that the sphere $S^3$ is the same as the unitary group $SU_2$ makes the presentation slightly more accessible  for the sphere than for the hyperboloid; for that reason the paper begins with the geometric preliminaries on $SU_2$ required for the  definition of  the symplectic form for the spherical Darboux curves. This material is presented in Section {\rm{2}}. The analogous material for the hyperbolic Darboux curves is presented in Section {\rm{3}}. 

 The second part of the paper begins with Section{\rm{4}} and is devoted to the Hamiltonian flow corresponding to $f_0=\half\int_0^L\kappa^2(s)\,ds$. This section also contains  a discussion of the Euclidean symplectic form and its connection to the existing results in the literature. Section {\rm{5}} deals with elastic curves and the soliton solutions for the non-linear Schroedinger's equation. The final section (Section {\rm{6}}) contains a brief discussion of the hierarchy of Poisson commuting functions and their connections to the hierarchies presented in ([14] and [18]).

\section{Darboux curves on the sphere and their symplectic form}

\subsection{Notations and geometric preliminaries}

  The three dimensional sphere $S^3=\{x\in \IR^4:x_0^2+x_1^2+x_3^2+x_4^2=1\}$  can be represented either by  the unit quaternions or by  the matrices $X=\bmt z&w\\-\bar{w}&\bar{z}\emt$ with $z=x_0+ix_1$, $w=x_2+ix_3$ and $|z|^2+|w|^2=1$. The most direct way to these representations is through the identification of points $x=(x_o,x_1,x_2,x_3)$ in $\IR^4$ with the  matrices
 \[X=x_0I+x_1E_1+x_2E_2+x_3E_3=\bmt x_0+ix_1 & x_2+ix_3\\-x_2+ix_3 & x_0-ix_1\emt\]  where \[E_1=\bmt i & 0\\0 & -i\emt,\,E_2=\bmt0 & 1\\-1 & 0\emt ,\,E_3=\bmt 0 & i\\i & 0\emt\]. Matrices $E_1,E_2,E_3$  conform to the relations
\[E_1^2=E_2^2=E_3^2=-I,\]
\[ E_1E_2=E_3,E_3E_1=E_2,E_2E_3=E_1\] and can be identified with $i,j,k$, the standard basis  for the imaginary quaternions.
The sphere $S^3$  then is equal to the  set of matrices $X=x_0I+x_1E_1+x_2E_2+x_3E_3$ whose  determinant $\Det(X)=1$. Alternatively, the sphere could also be  defined as the group of $2\times 2$ matrices $g$ with complex entries whose inverses $g^{-1}$ are equal to their Hermitian transposes $g^*$. This group is called the special unitary group and is traditionally denoted by $SU_2$. 
  
However, in this paper the special unitary group shall be denoted by $K$ for additional simplicity in notation. 
The Lie algebra of $K$ denoted  by $\fk$ consists of the matrices \[A=\bmt ix_1 & u_2+ix_3\\-x_2+ix_3 & -ix_1\emt\] with  $x_1,x_2,x_3$ arbitrary real numbers. The Lie bracket $[A,B]$ is defined by $[A,B]=BA-AB$ for $A$ and $B$ in $\fk$ . The tangent vectors at any point $X$ of $K$ shall be represented by matrices $XA$ with $A\in \fk$, and the tangent bundle $TK$ of $K$ shall be represented by the product $K\times \fk$ in terms of the left-invariant vector fields $X\rightarrow XA,\,A\in\fk$. 

{\definition Matrices 
\eq A_1={1\over 2}\bmt i&0\\0&-i\emt ,\,A_2=\bmt0&1\\-1&0\emt,\,A_3=\bmt0&i\\i & 0\emt
\ee
shall be known as the Pauli matrices.}

 The Pauli matrices form a basis for $\fk$ and  conform to the following Lie bracket relations:
\eq \label{paulibrackets}
 [A_1,A_2]=-A_3,\,[A_1,A_3=A_2,\,[A_2,A_3]=-A_1\ee

{\definition  The trace form is a quadratic form on $\fk$ defined by\eq
\langle A,B\rangle=-2Tr(AB)\ee
for any matrices $A,B$ in $\fk$, where $Tr(AB)$ denotes the trace of $AB$.}

It follows that\[\langle A,B\rangle=\sum u_iv_i\] whenever \[A=\sum u_iA_i \mbox{ and }B=\sum v_iA_i.\] Therefore, the trace form is positive definite on $\fk$, and the matrices $A_1,A_2,A_3$ form an orthonormal  basis in $\fk$. Then, \[||A||=\sqrt{\langle A,A\rangle}=\sqrt{u_1^2+u_2^2+u_3^2}\] for any matrix $A=\sum u_iA_i$.

 The trace form, in addition,  satisfies these invariance properties:\eq
\langle A,[B,C]\rangle=\langle[A,B],C\rangle\mbox{, and }\langle gAg^*,gBg^*\rangle=\langle A,B\rangle\ee
for any matrices $A,B,C$ in $\fk$, and any $g$ in $K$.  

The Riemannian structure of the sphere, inherited from the Euclidean inner product $(x,y)$ in $\IR^4$ is related to the trace form according to the following lemma.

{\lemma\label{riemmetric} Suppose that $\frac{dX_1}{ds}(0)=XA$ and $\frac{dX_2}{ds}(0)=XB$
denote tangent vectors  at $X$ represented by the curves $X_1(s)$ and $X_2(s)$ in $K$.  
 Then,\eq
(\frac{dx_1}{ds}(0),\frac{dx_2}{ds}(0))=4\langle A,B\rangle.\ee}

\prf For each matrix $X=x_0I+x_1E_1+x_2E_2+x_3E_3$ the  companion matrix $X^{\dagger}$ is defined by \[X^{\dagger}=x_0I-x_1E_1-x_2E_2-x_3E_3=\bmt x_0-ix_1 & -x_2-ix_3\\x_2-ix_3 & x_0+ix_1\emt\] Evidently, $X^{\dagger}=-X$ when $X$ belongs to $\fk$, and 
\eq\label{euctrace} \half (XY^{\dagger}+YX^{\dagger})=(x,y)I,\ee
as can be verified by an easy calculation.

Therefore,
\bean (\frac{dx_1}{ds}(0),\frac{dx_2}{ds}(0))I&=&\half((\frac{dX_1}{ds}(0))(\frac{dX_2}{ds}(0))^{\dagger}+(\frac{dX_2}{ds}(0))(\frac{dX_1}{ds}(0))^{\dagger})\\&=&
\half((XA)(XB)^{\dagger}+(XB)(XA)^{\dagger})=\half(AB^{\dagger}+BA^{\dagger})\\&=&
-\half(AB+BA)=4\langle A,B\rangle
\eean
\qed  

{\definition The group $K$ shall be considered a Riemannian manifold with its Riemannian metric
given by the trace form.}

It follows that the length $l(X)$ of any curve $X(s),\,s\in [0,L]$ is given by \[l(X)=\int_0^L||\Lambda (s)||\,ds ,\] where $\Lambda (s)=X(s)^*\frac{dX}{ds}(s)$. The length $l(X)$ is equal to half the length of the Euclidean metric in $\IR^4$ given by the usual formula $\int_0^L(\frac{dx}{ds}(s),\frac{dx}{ds}(s))^{\half}\,ds$

{\definition\label{covder0}
 The covariant derivative  of a curve of 
tangent vectors $v(s))=X(s)B(s)$ along  a curve 
$X(s)$ is defined by
\[
\frac{D_X}{ds}(V)(s)=X(s)(\frac{dB}{ds}(s)+{1\over 2}[A(s),B(s)])\]
where $A(s)=X^*(s)\frac{dX}{ds}(s)$.}

The reader can easily verify that the covariant derivative is equal to the orthogonal projection of the ordinary derivative in $\IR^4$ onto the tangent space od the sphere when the sphere is considered a subset of $\IR^4$. 

{\lemma\label{mixcovder} Suppose that $X(s,t)$ is a field of curves in $K$ with the infinitesimal directions \[A(s,t)=X^*(s,t)\frac{\partial X}{\partial s}(s,t)\mbox { and }B(s,t)=X^*(s,t)\frac{\partial X}{\partial t}(s,t)\]
Then,\eq
\frac {\partial A}{\partial t}(s,t)-\frac {\partial B}{\partial s}(s,t)+[A(s,t),B(s,t)]=0
\ee}

\prf  As in any Riemannian manifold,\[\frac{D_X}{ds}(\frac{\partial X}{\partial t})(s,t)=\frac{D_X}{dt}(\frac{\partial X}{\partial s})(s,t)\]
Therefore,
\[ 
X(s,t)(\frac{\partial B}{\partial s}(s,t)+{1\over 2}[B(s,t),A(s,t)])=X(s,t)(\frac{\partial A}{\partial t}(s,t)+{1\over 2}[A(s,t),B(s,t)])
\]
which implies that\[\frac {\partial A}{\partial t}(s,t)-\frac {\partial B}{\partial s}(s,t)+[A(s,t),B(s,t)]=0\]
\qed

 This equation is known as the zero-curvature equation ([5]).

\subsection {The orthonormal frame bundle and the Darboux curves}

The orthonormal frame bundle of $K$ is the totality of pairs $(X,F)$ where $X$ is a point in $K$ and where $F=(v_1,v_2,v_3)$ is an orthonormal frame at $X$.  Curves in the orthonormal frame bundle of $K$ shall be called  framed curves.  Basic to this paper is a set of framed curves defined by the following properties:

(a) Each curve is defined over a fixed interval [0,L], and conforms to the fixed initial condition \[X(0)=I,F(0)=(v_1(0),v_2(0),v_3(0))=(A_1,A_2,A_3).\]

(b) For each curve $(X(s),F(s))$ the orthonormal frame $F(s)=(v_1(s),v_2(s),v_3(s))$  defined  along $X(s)$ is adapted to  the curve $X(s)$ by the requirement that $v_1(s)=\frac{dX}{ds}$ for all $s\in [0,L]$.

{\definition Framed curves that satisfy conditions (a)
 and (b) above shall be called anchored Darboux curves.}

 For each Darboux curve $( X(s),F(s))$ condition (b) implies that $||\frac{dX}{ds}(s)||=1$, which then implies that $L$ is the length of the projected curve curve $X$.

Darboux curves  can be embedded in the space of curves in $K\times K$, the universal cover of the orthonormal frame bundle of $K$ through the following realization.
Each pair $(p,q)$ in $K\times K$ defines an orthonormal frame $(v_1,v_2,v_3)$ at $X=pq^*$  where \[v_1=pA_1q^*=pq^*(qA_1q^*),v_2=pA_2q^*=pq^*(qA_2q^*),v_3=pA_3q^*=pq^*(qA_3q^*)\] 

 Conversely, every orthonormal frame at a point $X\in K$  can be represented by the tangent vectors $v_1=XB_1,v_2=XB_2,v_3=XB_3$ for some matrices $B_1,B_2,B_3$ in $\fk$ that are orthonormal relative to the trace form. Then, there is matrix $q\in K$ such that

\[v_1=qA_1q^*,v_1=qA_2q^*,v_1=qA_3q^*\]
There are exactly two matrices $\pm q$ which satisfy the preceding equalities. Having found $q$, $p$ is uniquely defined by $p=Xq$.  

The fact that $K\times K$ is a double cover of the orthonormal frame bundle of the sphere does not matter for  the  subsequent exposition since  Darboux curves are  anchored at the identity. 

It follows from above that anchored Darboux curves  can be represented by curves $(p(s),q(s))$ in $K\times K$  that are  the solutions of 
\eq
\frac{dp}{ds}=p(s)P(s),\,\frac{dq}{ds}=q(s)Q(s)
\ee
 subject to the conditions that $p(0)=I$, $q(0)=I$, and that $P(s)-Q(s)=q(s)A_1q^*(s)$. 
The last condition reflects the fact that the curve of frames \[(X(s)q(s)A_1q^*(s),X(s)q(s)A_2q^*(s),X(s)q(s)A_3q^*(s))\] along the curve $X(s)$ in $K$  is adapted to the tangent vector $X(s)A(s)$ of $X(s)$
 through   the relation $X(s)A(s)=X(s)q(s)A_1q^*(s)$, which implies 
that $A(s)=q(s)A_1q^*(s)$.
   
Therefore anchored Darboux curves can be redefined as pairs of curves $(X(s),R(s))$ in $K\times K$ that are the solutions of  the initial value problem
 \eq \label{dar1}\frac{dX}{ds}(s)=X(s)\Lambda (s),\,\frac{dR}{ds}(s)=R(s)A(s)\ee,
\eq \label{dar2}X(0)=R(0)=I,\mbox { subject to }\Lambda (s)=R(s)A_1R^*(s)\ee.

{\definition\label{horizdarboux} The projection $X(s)$ of an anchored Darboux curve shall be called horizontal Darboux curve.}
 
    The horizontal Darboux curves are parametrized by three  arbitrary functions $u_1(s),u_2(s),u_3(s)$ (given by $A(s)=\sum u_i(s)A_i$ in (\ref{dar2})). These functions are related to the curvature $\kappa (s)$ and the torsion $\tau (s)$  of the projected curve $X(s)$ through the following equations:
\[\kappa (s)=u_3(s)\cos{\theta (s)}-u_2(s)\sin{\theta (s)},\,\frac{d\theta}{ds}(s)+u_1(s)=\tau (s)\] subject to the constraint \[u_2(s)\cos{\theta (s)}+u_3(s)\sin{\theta (s)}=0\]
 The angle $\theta$ represents the rotation of the Darboux frame relative to the Serret-Frenet frame in the plane perpendicular to the tangent vector of the projected curve. 

Recall that the Serret-Frenet frame $(v_1(s),v_2(s),v_3(s))$, generated by a curve $X(s)$, is defined by the following relations: $\frac{dX}{ds}=v_1$ and
\eq\label{serfr}
\frac{D_X}{ds}(v_1)=\kappa v_2,\frac{D_X}{ds}(v_2)=-\kappa v_1+\tau v_3,\frac{D_X}{ds}(v_3)=-\tau v_2
\ee

 In the representation of the frame bundle as $K\times K$, the Serret-Frenet frame is  represented by the matrix $R(s)$ through the relations\[v_1(s)=X(s)(R(s)A_1R^*(s)),v_2(s)=X(s)(R(s)A_2R^*(s)),v_3(s)=X(s)(R(s)A_3R^*(s))\]
It now follows from Definition \ref{covder0} that \bean \frac{D_X}{ds}(v_1)&=& X(R[A_1,A]R^*)=-u_3 X(RA_3R^*)+u_2X(RA_2R^*)\\&=& u_3v_2-u_3v_3\eean
  from which it follows that $u_3=\kappa $ and $u_2=0$. Furthemore,
\bean
\frac{D_X}{ds}(v_2)(s)&=& X(s)R(s)([A_2,A(s)]+{1\over 2}[A_1,A_2])R^*(s))\\&=&=-u_3(s)X(s)(R(s)A_1R^*(s))+(u_1(s)-{1\over 2})X(s)(R(s)A_3R^*(s))\\&=&=-u_3(s)v_1(s)+(u_1(s)-{1\over 2})v_3(s)\eean
implies that $\tau (s)=u_1(s)-{1\over 2}$.

Therefore the frame curve $R(s)$ corresponds to the Serret-Frenet frame  whenever $u_2(s)=0$, in which case $u_3(s)=\kappa (s)$ and $u_1(s)=\tau (s)+{1\over 2}$.

More generally,

{\theorem\label{curvdarb} (a) For any solution $(X(s),R(s))$ of equations (\ref{dar1}) and (\ref{dar2}) with $A(s)=\sum u_i (s)A_i$,
\[\kappa^2(s)=u_2^2(s)+u_3^2(s)\]
where $\kappa (s)$ denotes the curvature of the base curve $X(s)$.

(b) Every curve $X(s)$ which satisfies $||\frac{dX}{ds}(s)||=1$,  $X(0)=I,$ and $\frac{dX}{ds}(0)=A_1$ is the projection of an anchored Darboux curve.}

\prf \bean\kappa ^2(s)&=&||\frac{D_X}{ds}(\frac{dX}{ds})(s)||^2\\
&=&||X(s)(\frac{d}{ds}(R(s)A_1R^*(s))||^2=||R(s)[A_1,A(s)]R^*(s)||^2\\
&=&||-u_2(s)A_3+u_3(s)A_2||^2 = u_2^2(s)+u_3^2(s)
\eean
To prove part (b) let $\Lambda (s)=X^*(s)\frac{dX}{ds}(s)$. It follows from Lemma \ref{riemmetric} that $||\Lambda (s)||=1$. Since $K$ acts transitively on $S^3$ it follows that there exists a curve $R(s)$ such that $R(s)A_1R^*(s)=\Lambda (s)$. Then $(X(s), R(s))$ is the desired  anchored Darboux curve.
\qed

{\bf{Remark}}. The reader may recall that not all curves can be lifted to the Serret-Frenet frames.  Darboux frames do not have that disadvantage, as the preceding theorem shows, but  the uniqueness is lost. In fact, for any Darboux curve $(X(s),R(s))$ $(X(s),R(s)\phi (s))$ is another anchored Darboux curve that projects onto $X(s)$ provided that $\phi (s)A_1\phi ^* (s)=A_1$. 

{\definition\label{darbs} The set of anchored Darboux curves will be denoted by $\mathcal{D}_s(L)$, while $Horiz(\mathcal{D}_s)(L)$ will denote the set of all horizontal Darboux curves.}

{\definition\label{perdarb} An anchored Darboux curve $(X(s),R(s))$ will be called frame- periodic if $R(0)=R(L)$, i.e., if $R(s)$ is a closed curve. The space of frame-periodic Darboux curves will be denoted by $\mathcal{PD}_s(L)$. Then, $Horiz(\mathcal{PD}_s)(L)$ will denote the set of all horizontal frame-periodic Darboux curves.} 

If $(X(s),R(s))$ is frame-periodic then both matrices $\Lambda (s)$ and $A(s)$ defined by equations (\ref{dar1}) and (\ref{dar2}) are smoothly periodic, because closed solutions of differential equations with smooth data are smoothly periodic. However, the projected curve $X(s)$ need not be closed. 

On the other hand,  the Serret-Frenet frame of any smoothly periodic curve $X(s)$ is periodic, and therefore, all smoothly periodic curves in $K$ are the projections of frame-periodic Darboux curves.

\subsection{ Darboux curves as Fr\'echet manifolds}

 On the basis of the general theory developed in ([6]) both  $\mathcal{D}_s(L)$ and $\mathcal{PD}_s(L)$  may be considered as infinite-dimensional Fr\'echet manifolds. Recall that
a  topological  Hausdorff vector space $V$ is called a Fr\'echet space if its topology is induced by a countable family of semi-norms $p_n$, and if it is complete relative to 
 the semi-norms in $\{p_n\}$. A Fr\'echet manifold is defined as follows:

{\definition A Fr\'echet manifold  is a topological
Hausdorff space equipped with an atlas whose
charts take values in open subsets of a Fr\'echet
space $V$  such that any change of coordinate 
charts is smooth.}

The paper of R.S. Hamilton ([6]) singles out an important class of Fr\'echet manifolds, called tame, in which  the implicit function theorem is true. One of the main theorems in ([6]) is that 
  the set of smooth mappings from a compact manifold interval into a  finite-dimensional Riemannian manifold $M$ is a  tame Fr\'echet manifold. It therefore follows from the implicit function theorem that closed subsets of tame Fr\'echet manifolds ${\cal M}$, defined by  the zero sets of finitely many smooth functions on ${\cal M}$ are  tame sub-manifolds of ${\cal M}$. Since both $\mathcal{D}_s(L)$ and $\mathcal{PD}_s(L)$  are particular cases of the above situation, it follows that they are both tame Fr\'echet manifolds and the same applies to their horizontal projections $Horiz(\mathcal{D}_s)(L)$ and $Horiz(\mathcal{PD}_s)(L)$. 

Tangent vectors and tangent bundles of Fr\'echet manifolds are defined in the same manner as for finite dimensional manifolds. In particular tangent vectors at a point $x$ in  a Fr\'echet manifold ${\cal M}$  are the equivalence classes
of curves $\sigma(t)$ in ${\cal M}$ all emanating from $x$ (i.e.,$\sigma(0)=x$), and
 all having the same tangent
vector ${d\sigma\over dt}(0)$ in each equivalence class. The set of all tangent vectors at $x$ denoted by $T_x{\cal M}$ constitutes the tangent space at $x$. 

  The tangent bundle of a Fr\'echet manifold ${\cal M}$ is a Fr\'echet manifold. A vector field
$X$ on ${\cal M}$ is a smooth mapping from $\cal M$
into the tangent bundle $T\cal M$  such that $X(s)\in
T_x\cal M$ for each $x\in\cal M$. On tame  Fr\'echer manifolds vector fields can be defined as derivations in the space of smooth functions on ${\cal M}$.

In the formalism  of Fr\'echet manifolds tangent vectors in each of $Horiz(\mathcal{D}_s)(L)$ and $Horiz(\mathcal{PD}_s)(L)$ are given by the following theorem.

{\theorem \label{tgt} (a). The tangent space ${T_X(Horiz(\cal D}_s )(L)$  at an anchored  horizontal Darboux curve
$(X(s)$ consists of curves $v(s)=X(s)V(s)$ with $V(s)$ the solution of
\begin{equation}\label{variat}
\frac{dV}{ds}(s)=[\Lambda (s),V(s)]+U(s)
\end{equation}
 such that $V(0)=0$, where $U(s)$ is a
curve in $\fk$ subject to  the conditions that $U(0)=0$ and $\langle \Lambda (s),U(s)\rangle=0$. The matrix $\Lambda (s)$ is the tangent vector of $X$, i.e.,
 \[\frac{dX}{ds}(s)=X(s)\Lambda (s)\]

(b). $v(s)=X(s)V(s)$  is a tangent vector at an anchored frame-periodic horizontal Darboux curve $X(s)$ if in addition to the above, the curve $U(s)$ is smoothly periodic having the period equal to $L$.}

\prf Let $Y(s,t)$ denote a family of  anchored  horizontal Darboux curves such that $Y(s,0)=X(s)$. Then, $v(s)=\frac{\partial Y}{\partial t}(s,t)_{t=0}$ is a tangent vector at $X(s)$ for which $v(0)=0$ since the curves $Y(s,t)$ are anchored.

Let $Z(s,t)$ and $W(s,t)$ denote the matrices defined by 
\[Z(s,t)=Y(s,t)^*\frac{\partial Y}{\partial s}(s,t)\mbox {, }W(s,t)=Y(s,t)^*\frac{\partial Y}{\partial t}(s,t)\]
It follows that $\Lambda (s)=Z(s,0)$, $V(s)=W(s,0)$ and that $V(0)=0$. Furthemore,\[
\frac{\partial Z}{\partial t}(s,t)-\frac{\partial W}{\partial s}(s,t)+[Z(s,t),W(s,t)]=0\]

as a consequence of Lemma \ref{mixcovder}.  The above equation reduces to
\[
\frac{dV}{ds}(s)=[\Lambda (s),V(s)]+U(s)
\]

 when $t=0$ provided that $U(s)=\frac{\partial W}{\partial s}(s,0)$.

Since the curves $s\rightarrow Y(s,t)$ are Darboux for each $t$,\[\langle Z(s,t),Z(s,t)\rangle=1\mbox{, and } Z(0,t)=A_1\]
Therefore,
\[\langle \frac{\partial Z}{\partial t}(s,t),Z(s,t)\rangle=0\mbox{, and } \frac{\partial Z}{\partial t}(0,t)=0\]
which implies that $\langle \Lambda (s),U(s)\rangle=0$ and $U(0)=0$. 

It remains to show that any curve $V(s)$ in $\fk$ that satisfies (\ref{variat}) can be realized by the perturbations $Y(s,t)$ used above. So assume that $V(s)$ be any solution of (\ref{variat}) generated by  a curve $U(s)$ with $U(0)=0$ that satisfies $\langle \Lambda (s),U(s)\rangle=0$.

Let $\phi (t)$ denote any smooth function such that $\phi (0)=0$ and $\frac{d\phi}{dt}(0)=1$. Define
\[
Z(s,t)={1\over {1+\phi ^2 (t)\langle U(s),U(s)\rangle }}(\Lambda (s)+\phi (t)U(s))
\]
Evidently $Z(0,t)=A_1$ for all $t$, and an easy calculation shows that \newline$\langle Z(s,t),Z(s,t)\rangle =1$. Therefore $Y(s,t)$, the solution of \[\frac{\partial Y}{\partial s}(s,t)=Y(s,t)Z(s,t)\] with $Y(0,t)=I$ corresponds  to an anchored  horizontal Darboux curve for each $t$. Since $U(s)=\frac{\partial Z}{\partial t}(s,0)$ our proof of part (a) is finished.

To prove part (b) we shall assume that the curves $Y(s,t)$ used in part (a) belong to $Horiz({\cal {PD}}_s )(L)$. Then, curves $s\rightarrow Z(s,t)$ are L-periodic for each $t$, and therefore, $U(s)=\frac{\partial Z}{\partial t}(s,0)$ is periodic with period $L$.

\qed

\subsection{ The symplectic structure of horizontal Darboux curves}

The basic notions of symplectic geometry  of infinite-dimensional Fr\'echet manifolds are  defined through differential forms in the same manner as for the finite-dimensional situations.
In particular, differential forms $\omega$ of degree n are mappings
\[ \omega :\underbrace{{\cal X(M)}\times \cdots \times 
{\cal X(M)}}_{n}\to C^{\infty}({\cal M}) \]
that are $C^{\infty}({\cal M})$ multilinear and
skew-symmetric. Here ${\cal X(M)}$ denotes the space of all smooth vector fields on $\cal{M}$.

\begin{definition}\label{extder} The exterior derivative $d\omega$ of a form of degree $n$ 
is a differential form of degree $n+1$ defined by
\begin{eqnarray*}
	d\omega (X_1,\dots, X_{n+1})& = &
	\sum\limits_{i=1}^{n+1}(-1)^{i+1}X_i
	(\omega (X_1,\dots,\hat X_1,\dots, X_n))\\
	&&-\sum\limits_{i<j}(-1)^{i+j} \omega 
	 ([X_i,X_j],\dots \hat{X_i},\dots \hat{X_j},
	X_{n+1})\nonumber
\end{eqnarray*}\end{definition}
where the roof sign above an entry indicates its
absence from the expression (i.e., $w(\hat{X_1},X_2)
=w(X_2)$ and $w(X_1,\hat{X_2})=w(X_1)$).

A differential form $\omega$ is said to be closed if
its exterior derivative $d\omega$ is equal to zero.

\begin{definition}\label{symp} A differential form $\omega$ of degree 2 is
said to be symplectic whenever
it is closed and non-degenerate, in the sense that  the induced
form $(i_{X}\omega )(Y)=\omega (X,Y)$ is non-zero for each non-zero vector field $X$.\end{definition}

The differential $df$   of a smooth function $f$ is a form  of degree 1 defined by $df(v)= {d\over dt}
f\circ \sigma(t)|_{t=0}$ for any smooth curve
in ${\cal M}$ such that $\sigma(0)=x$, and $\frac{d\sigma }{dt}(0)=v$.

In finite dimensional symplectic manifolds with a symplectic form $\omega$ there
is a unique vector field $X_f$ such that $df=i_{X_f}\omega $.
$X_f$ is called the Hamiltonian vector field induced by 
$f$, and $f$ is called the Hamiltonian of $X_f$. However, in 
infinite dimensional manifolds it may happen 
that the form $df$ is not equal to $i_{X}w$ for any
$X\in{\cal X(M)}$. This is due to the fact that the cotangent bundle of an
 infinite dimensional Fr\'echet space is never a Fr\'echet manifold.  Nevertheless,

\begin{definition}\label{hamvect} A vector field $X$  is said to be Hamiltonian if there exists a smooth 
function $f$ such that
\[ df(Y)=\omega (X,Y)\]
for all vector fields $Y$ on ${\cal M}$. The dependence of $X$ on $f$ shall be noted explicitly by $X_f$.\end{definition}

The manifold consisisting of  horizontal Darboux curves posesses a natural differential form $\omega$ defined as follows:

Let $v_1(s)=X(s)V_1(s)$ and $v_2(s)=X(s)V_2(s)$ denote any tangent vectors at a horizontal Darboux curve $X(s)$ that is defined by $\frac{dX}{ds}(s)=X(s)\Lambda (s)$.
 According to (\ref{variat}) there exist unique curves $U_1(s)$ and $U_2(s)$ such that 
 \[U_i(0)=0,\, \langle \Lambda (s),U_i(s)\rangle=0\] 
and 
\eq\label{variat1}U_i(s)=\frac{dV_i}{ds}(s)-[\Lambda (s),V_i(s)],\,i=1,2.\ee Then $\omega$ is given by 

\begin{equation}\label{dsymp}
\omega _{\Lambda}(V_1, V_2)=-\int^L_0{\langle \Lambda (s), [U_1(s),U_2(s)]
\rangle }\,ds
\end{equation} 
{\it {Remark.}} As in finite dimensional situations the choice of sign is a matter of convention. The justification for the above choice of sign will  be given later on in the paper.
 
{\theorem\label{sympthm} Both  $Horiz({\cal{D}})_s(L)$ and  $Horiz({\cal{PD}})_s(L)$
are symplectic Fr\'echet manifolds relative to  $\omega$ defined by (\ref{dsymp}).}

The following lemmas will be useful for the proof of the theorem.
{\lemma \label{triplebracket}\[[A,[B,C]]=\langle A,C\rangle B -\langle A,B\rangle C)\]
 for any elements $A,B,C$ in $\fk$.}

We leave the proof to the reader.

{\lemma \label{control}Suppose that $v(s)=g(s)V(s)$ is any tangent vector at a horizontal Darboux curve $g(s)$.
Let $U(s)$ be defined by (\ref{variat}). Then there exists a curve $C(s)$ in $\fk$ such that 
\[U(s)=[\Lambda (s),C(s)]\]
\prf  The mapping $C\rightarrow [\Lambda (s),C(s)]$ restricted to the orthogonal complement of $\Lambda (s)$ is surjective. Since $U(s)$ is orthogonal to $\Lambda (s)$ the proof follows.
\qed

\bf{Proof of the theorem. }}  The proof is the same for each of $Horiz({\cal{D}})_s(L)$ and  $Horiz({\cal{PD}})_s(L)$  and will be presented formally without any reference to the underlying space.

Evidently $\omega $ is skew-symmetric. To show that it is non-degenerate, assume that 
for some tangent vector $v_1=gV_1$ at $g$,  $\omega _{\Lambda} (V_1, V)=0$ for all tangent vectors 
$v=gV$ at $g$. Let $U(s)$ and $U_1(s)$ correspond to $V(s)$ and $V_1(s)$ as in (\ref{variat1}).
Then take $U(s)=[\Lambda (s),U_1(s)]$. Evidently $ U(0)=0$, and $\langle \Lambda (s),U(s)\rangle=0$, and therefore the corresponding vector $v=gV$ belongs to the tangent space at $g$.  It follows from
Lemma \ref{triplebracket} that\[ [U_1[\Lambda,U_1]=\langle U_1,U_1\rangle\Lambda =||U_1||^2\Lambda\]
Therefore,
\[\langle\Lambda (s),[U_1(s),U(s)]\rangle=||\Lambda (s)||^2||U_1(s)||^2=||U_1(s)||^2\]
which implies that $U_1(s)=0$ for all $s$, since $0=\omega_{\Lambda}(V_1,V)=\int_0^L ||U_1(s)||^2\,ds$.

But then (\ref{variat1}) implies that $V_1(s)=0$. Hence, $\omega$ is non-degenerate.

To show that $\omega$ is closed let $v_i (s)=g(s)V_i (s)$,\, $1\leq i\leq 3$ denote any three tangent vectors at a  fixed Darboux curve $g(s)$. It is required to show (Definition \ref{extder}) that\eq\label{cyclic}
d\omega (X_1,X_2,X_3)= 
	\sum_{cyclic}X_i(\omega (X_j,X_k)+\sum_{cyclic}\omega([X_i,X_j],X_k)=0\ee
where $X_i$ denote any vector fields such that $X_i(g)=v_i$ for each $i=1,2,3$.

  To make use of the above formula it becomes necessary to identify vectors $v_i$ with vector fields $X_i$. For the purposes of this identification it will be convenient to  make a slight change in the notations and write $\Lambda_z(s)$ instead of $\Lambda (s)$ for the tangential direction $\Lambda(s)$ of the horizontal Darboux curve $z(s)$. 

Vector fields $X_i$ will be induced by the matrices $U_i=\frac{dV_i}{ds}(s)-[\Lambda ,V_i]$ through the following relations: for each horizontal Darboux curve $z(s)$ let $Z_i(s)$ denote the solution of \[\frac{dZ_i}{ds}(s)=[\Lambda (z)(s),Z_i(s)]+[\Lambda (z)(s),U_i(s)]\] with $Z_i(0)=0$ for each $i=1,2,3$.  Then  \[X_i(z)(s)=z(s)Z_i(s),\,i=1,2,3\] 

To show that (\ref{cyclic}) is valid let $C_i(s)$ denote the curves such that $U_i(s)=[\Lambda (s),C_i(s)]$ as in Lemma \ref{control}. Then, $\frac{dV_i}{ds}=[\Lambda ,V_i+C_i]$ and  an easy calculation based on Jacobi's identity yields\[\frac{d}{ds}([V_i,V_j])=[\Lambda ,[V_i,V_j]+[[\Lambda,C_i],V_j]+[V_i,[\Lambda,C_j]]\]
Therefore,
\bean
\sum_{cyclic}\omega([X_i,X_j],X_k)&=&\sum_{cyclic}\int_0^L\langle\Lambda,[[[\Lambda,C_i],V_j]+[V_i,[\Lambda,C_j]],[\Lambda,C_k]]\,ds\\
&=&\sum_{cyclic}\int_0^L\langle\Lambda,[(\langle V_i,\Lambda\rangle C_j-\langle V_j,\Lambda\rangle C_i)+(\langle V_j,C_i\rangle \Lambda -\langle V_i,C_j\rangle \Lambda ),[\Lambda,C_k]]\,ds\\&=&\sum_{cyclic}\int_0^L\langle V_j,\Lambda\rangle\langle C_i,C_k\rangle -\langle V_i,\Lambda\rangle\langle C_j,C_k\rangle\,ds=0
\eean

The calculations involving $X_i(\omega (X_j,X_k)$  in (\ref{cyclic}) require additional notations. Let $t\rightarrow z_i(s,t)$ denote the integral curves of the vector field $X_i$ that originate at $g(s)$ for $t=0$, and let \[
\frac{\partial z_i}{\partial t}(s,t)=z_i(s,t)Z_i((z_i(s,t))\mbox{ and }\frac{\partial z_i}{\partial s}(s,t)=z_i(s,t)\Lambda_i (z_i(s,t)).\]
For simplicity of notation let $Z_i(z_i(s,t))$ and $\Lambda _i(z_i(s,t))$ be denoted by $Z_i(s,t)$ and $\Lambda _i(s,t)$ . From Lemma \ref{mixcovder}\eq\label{flow}\frac{\partial \Lambda _i}{\partial t}-\frac{\partial Z_i}{\partial s}+[\Lambda_i ,Z_i]=0\ee
 
 which at $t=0$ reduce to \[U_i-\frac{dV_i}{ds}+[\Lambda ,V_i]=0\]
 As in the preceeding calculation $U_i$ will be represented by $U_i=[\Lambda ,C_i]$.

Then, 
\bean
X_i(\omega (X_j,X_k)&=&\frac{\partial }{\partial t}\int_0^L\langle \Lambda_i (s,t),[[\Lambda_j (s,t),C_j],[\Lambda_k (s,t),C_k]] {\,}ds|_{t=0}\\
&=&\int_0^L\langle \frac{\partial\Lambda_i}{\partial t} (s,t),[[\Lambda_j (s,t),C_j],[\Lambda_k (s,t),C_k]]\rangle{ \,}ds|_{t=0}\\
&+&\int_0^L\langle \Lambda_i (s,t),[[\frac{\partial\Lambda_j}{\partial t} (s,t),C_j],[\Lambda_k (s,t),C_k]]\rangle{ \,}ds|_{t=0}\\
&+&\int_0^L\langle \Lambda_i (s,t),[[\Lambda_j (s,t),C_j],[\frac{\partial \Lambda_k}{\partial t} (s,t),C_k]]\rangle{ \,}ds|_{t=0}\\
&=&\int_0^L\langle U_i(s),[U_j(s),U_k(s)]\rangle\,ds\\& +&\int_0^L \langle \Lambda (s),(
 [[U_j(s),C_j],[\Lambda (s),C_k]]+[[\Lambda (s),C_j],[U_k(s),C_k]]\rangle){\,}ds
\eean

Since $\langle U_i(s),[U_j(s),U_k(s)]\rangle$  is equal to the volume of the parallelopiped with sides $ U_i(s),U_j(s),U_k(s)$, \[\langle U_i(s),[U_j(s),U_k(s)]\rangle=0\] because each side  of the parallelopiped is in the plane orthogonal to $\Lambda (s)$. Therefore,\[\int_0^L\langle U_i(s),[U_j(s),U_k(s)]\rangle\,ds=0\]
 
It remains to show that\[\sum_{cyclic}\int_0^L \langle \Lambda (s),(
 [[U_j(s),C_j],[\Lambda (s),C_k]]+[[\Lambda (s),C_j],[U_k(s),C_k]]\rangle){\,}ds=0\]

It follows that
 
 \[
 \langle \Lambda (s),(
 [[U_j(s),C_j],[\Lambda (s),C_k]]+[[\Lambda (s),C_j],[U_k(s),C_k]]\rangle)
=\]\[\langle \Lambda (s),(\langle [U_j(s),C_j],\Lambda (s)\rangle C_k-\langle [U_j(s),C_j],C_k\rangle \Lambda (s)
+\langle [U_j(s),C_j],\Lambda (s)\rangle C_k-\langle [U_k(s),C_k],C_j\rangle \Lambda (s)\rangle )=\]
\[\langle[U_k(s),C_k(s)],C_j\rangle -\langle[U_j(s),C_j],C_k\rangle =0\]

because $ \langle[U_k(s),C_k(s)],C_j\rangle$ is the volume of the parallelopiped with sides $U_k(s),C_k(s),C_j$   all of which are  in the plane orthogonal to $\Lambda (s)$.
 Therefore $\omega $ is closed, and hence symplectic.

\qed

\section{ The symplectic structure of hyperbolic Darboux curves}

\subsection{ The hyperboloid and its frame bundle}

 Similar to the sphere, the hyperboloid can be represented in several ways. To make easy transitions from the sphere to the present situation it will be most convenient to represent the hyperboloid as the homogeneous manifold $SL_2(C)/SU_2$. The representation is done as follows.

Each vector in $x$ in $\IR^4$, considered as a column vector with coordinates $x_0,x_1,x_2,x_3$, can be represented  by a Hermitian matrix \[X=\bmt x_0+x_1 & x_2+ix_3\\x_2-ix_3 & x_0-x_1\emt\].
If $(a,b)_h$  denotes the Lorentzian quadratic form in $\IR^4$ defined by \eq\label{lorentz}(a,b)_h=a_0b_0-(a_1b_1+a_2b_2+a_3b_3)
\ee
then the Lorentzian unit sphere $(x,x)_h=1$ is a hyperboloid of two sheets and corresponds to the  space of Hermitian matrices $X$ whose determinant is equal to $1$.

 {\definition\label{hyper}The connected component  of the Lorentzian  unit sphere defined by $x_0 >0$ shall be  referred to as the hyperboloid and will be denoted by $\HI^3$. As a subset of the Hermitian matrices, the hyperboloid $\HI^3$ is equal to the space of positive-definite Hermitian matrices $X$ whose determinant is equal to $1$.}
 
 The orthonormal frame bundle of $\HI^3$ is equal to $SO(1,3)$, the  matrix group that leaves  the  Lorentzian form (\ref{lorentz})  invariant (it is understood  here that $SO(1,3)$ acts  on the points of $\IR^4$ by  the matrix multiplications). The  restriction of this action to $\HI^3$ is transitive and consequently, $\HI^3$ can be identified with the orbit  of $SO(1,3)$ through any of its points $\hat x$. In this identifications points of $\HI^3$ are represented by the first columns of the matrices $R\in SO(1,3)$ when $\hat x=e_0$. The other columns of the matrix are identified with an orthonormal frame $v_1,v_2,v_3$ at the base point $x=Re_0$. More precisely, the correspondence between matrices $R\in SO(1,3)$  and points of the othonormal frame bundle of $\HI^3$ is given by the following relations:\eq\label{frame}x=Re_0,\,v_1=Re_1,\,v_2=Re_2,\,v_3=Re_3\ee
The  frame $e_1,e_2,e_3$ at the identity induces an orientation on the space of frames and identifies $SO_0 (1,3)$, the connected component of $SO(1,3)$ that contains the identity, as the positively oriented orthonormal frame bundle of $\HI^3$.

 The subgroup $H$ of $SO_0(1,3)$ that fixes $e_0$ is isomorphic to $SO_3(R)$ and acts by the right multiplications on each fiber of the projection map $\pi (R)=Re_0$. In the language of the principal bundles the above construction identifies $SO(1,3)$ as the principal $H$ bundle with $\HI^3$ the base space under the projection map $\pi$.
 
 The transition from the spherical case to the hyperbolic case is more direct if instead of  the realization  $SO(1,3)/SO_3(R)$ the hyperboloid is realized as the quotient  $SL_2(C)/SU_2$. The passage from $SO(1,3)$ to $SL_2(C)$  is obtained through the identification of points $x=(x_0,x_1,x_2,x_3)$  in $\IR^4$ with Hermitian matrices  \[X=\left (\begin{array}{cc}x_0+x_1 & x_2+ix_3\\x_2-ix_3 & x_0-x_1\end{array}\right )\] and a homomorphism $\Phi (p)=R$ from $SL_2(C)$ onto $SO(1,3)$ defined by :

 \begin{equation}\label{hisom}
pXp^*=Y\mbox{ if and only if }Rx=y\mbox{, }x\in R^4
\end{equation}

 That the matrix $R$ defined by (\ref{hisom}) is in $SO(1,3)$ can be easily seen through the following lemma. 
{\lemma\label{hyptr} Suppose $X$ and $Y$ are Hermitian matrices  that correspond to points $x$ and $y$ in $\IR^4$. Then,

(a). ${1\over 2}(XY^\dagger +YX^\dagger )=(x,y)_h I$.

(b). $(gX)(gY)^\dagger +(gY)(gX)^\dagger =XY^\dagger +YX^\dagger $ for any $g$ in $SL_2(C)$.}
 
In the statement of the lemma $X^{\dagger}$ denotes the companion matrix  of a Hermitian matrix $X$. Recall that \[X^{\dagger}=\bmt x_0-x_1 & -x_2-ix_3\\-x_2+ix_3 & x_0+x_1\emt \mbox{ for any } X=\bmt x_0+x_1 & x_2+ix_3\\x_2-ix_3 & x_0-x_1\emt.\] 
 The proof of this lemma  will be left to the reader.
 
 The homomorphism $\Phi$ given by (\ref{hisom}) shows that $SL_2(C)$ is a double cover of $SO(1,3)$ since $ker(\Phi)=\{\pm I\}$. Consequently, the restriction of $\Phi$ to $H=\{1\}\times SO_3(R)$ proves $SU_2$ is a double cover of $SO_3(R)$.
 Additionally,
    the action $(g,X)\rightarrow gXg^*$ identifies $\HI^3$ as the quotient $SL_2(C)/SU_2$ via the orbit through the identity. 

The frame $v_1=Re_1,v_2=Re_2,v_3=Re_3$ at $x=Re_0$ corresponds to the matrices \eq\label{hypframe}V_1=pE_1p^*,\,V_2=pE_2p^*,\,V_3=pE_3p^*\ee  where
\[E_0=I,E_1=\bmt-1 & 0\\0 & 1\emt,E_2=\bmt0 & 1\\1 & 0\emt,E_3=\bmt 0 & i\\-i & 0\emt\]
are  the Hermitian matrices that correspond to the standard basis $e_0,e_1,e_2,e_3$ of $\IR^4$.

\subsection{Lie algebras and the trace form}

The Lie algebra of $SL_2(C)$ consisting of $2\times 2$ complex matrices of 
trace zero shall be denoted by $\fg$. The vector space of all 
Hermitian matrices in ${\fg}$
will be denoted by ${\fh}$ and the Lie subalgebra of all skew-Hermitian matrices in ${\fg}$ 
  will be denoted by ${\fk}$ in accordance with the notations used earlier in the paper. Then,
\[\fg=\fh+\fk\]

 and the following Lie algebraic relations hold
\begin{equation}\label{cartan}
[\fh,\fh]=\kappa,\quad [\fh, \fk]=\fh {\rm \quad and \quad}
[\kappa,\fk]=\fk
\end{equation}

The relation $[\fh,\fk]=
\fh$ easily implies that
$\vp A \vp^{-1} \in \fh$ for any $\vp \in K$
and any $A$ in $\fh$, which further implies that
the same holds for absolutely continuous curves
\begin{equation}
\vp(s) A(s) \vp^{-1}(s)
\end{equation}
with $\vp(s)$ a curve in $K$ and $A(s)$ a curve in $\fh$.

{\definition The matrices
\[B_1={1\over 2}\left (\begin{array}{cc}1 & 0\\0 & -1\end{array}\right ),
B_2={1\over 2}\bmt 0 & -i\\i & 0\emt,
B_3={1\over 2}\left (\begin{array}{cc}0 & 1\\1 & 0\end{array}\right ) 
\] shall be called the Hermitian Pauli matrices.}

 The Hermitian Pauli matrices are related to
the skew-Hermitian matrices $A_1,A_2,A_3$ defined earlier by the following simple correspondences \eq\label{hermpauli}A_1=iB_1,A_2=iB_2,A_3=iB_3\ee
Matrices $B_1,B_2,B_3,A_1,A_2,A_3$  form a basis for $\fg$, and
the reader can readily verify that their Lie brackets conform to the following Lie bracket table:

\centerline{\begin{tabular}{|c|c|c|c||c|c|c||}
\cline{2-7}
\omit [ , ] & \multicolumn{1}{|c|}{$A_1$} 
	& $A_2$ & $A_3$ & $B_1$ & $B_2$ & $B_3$ \\ \hline
%\omit\hfill[ , ] \hfill\vline& $A_1$ & $A_2$ & $A_3$ 
	%& $B_1$ & $B_2$ & $B_3$ \\ \hline
$A_1$ & $0$ & $-A_3$ & $A_2$ & $0$ & $-B_3$ & $B_2$ \\ \hline
$A_2$ & $A_3$ & $0$ & $-A_1$ & $B_3$ & $0$ & $-B_1$ \\ \hline
$A_3$ & $-A_2$ & $A_1$ & $0$ & $-B_2$ & $B_1$ & $0$ \\ \hline \hline
$B_1$ & $0$ & $-B_3$ & $B_2$ & $0$ & $A_3$ & $-A_2$ \\ \hline
$B_2$ & $B_3$ & $0$ & $-B_1$ & $-A_3$ & $0$ & $A_1$ \\ \hline
$B_3$ & $-B_2$ & $B_1$ & $0$ & $A_2$ & $-A_1$ & $0$ \\ \hline
\end{tabular}}

\centerline{Table 1}
It may
 be helpful for some of the subsequent calculations to note
the following relations:
{\lemma \label{vectprod}(a)
	If $A$ and $B$ are skew-Hermitian  with  $A=\sum\limits_{i=1}^3a_iA_i$ and $B=
	\sum\limits_{i=1}^3b_iA_i$, then $[A,B]=\sum\limits_{
	i=1}^3c_iA_i$ where $ c$ is the vector product $b\times a$ in $R^3$.
	
(b) If  $A$ is skew-Hermitian with $A=\sum \limits_{i=1}^3a_iA_i$, and $B$ is Hermitian with $B=\sum \limits_{i=1}^3b_iB_i$,
	then again $[A,B]=\sum\limits_{i=1}^3 c_iB_i$ with $c=b\times a$.

	(c) However,  if both $A$ and $B$ are Hermitian with $A=\sum \limits_{i=1}^3a_iB_i$ and $B=\sum\limits_{i=1}^3b_iB_i$ then
	$[A,B]=\sum \limits_{i=1}^3c_iA_i$ with $c=a\times b$. Thus for Hermitian matrices
 the order in the cross product is reversed.}

\begin{definition}   For any matrices $A,B $ in $\fg$, $\leangle A, B \riangle$ will denote
  the trace of  $2(AB)$. We shall refer to this quadratic form as the trace form on $sl_2(C)$.\end{definition}
 
 The Hermitian Pauli matrices $B_1,B_2,B_3$ form an orthonormal basis in $\fh$ relative to the trace form, which together with (\ref{hermpauli})  implies that\[\langle A_i,A_j\rangle =-\delta_{ij}\mbox{ and } \langle A_i,B_j\rangle =i\delta_{ij}\] 

{\bf{Remark.}} The above relations imply that the  restriction of the trace form to the matrices in $\fk$ is the negative of the trace form used in the first part of the paper. For the most part of the paper this ambiguity in terminology will not matter since the meaning will be clear from the context, but in the instances that require clarification the meaning will be made explicit.
 
It follows from above that \[\langle A,A\rangle =-(a_1^2 +a_2^2 +a_3^2),\,
\langle B,B\rangle =b_1^2+b_2^2+b_3^2,\,\mbox { and }
\langle A,B\rangle =i(a_1b_1+a_2 b_3+a_3 b_2)\]
for any  matrices $A=\sum a_jA_j$ and $B=\sum b_jB_j$. The above  imply that
\[\fh \cong \fk^*\]
 where $\fk^*$ denotes the dual of $\fk$, a fact of central importance for
 the results that follow.

\subsection{Hyperbolic Darboux curves}

 As in the first part of the paper  horizontal Darboux curves shall are naturally introduced through the language of the principal bundles. Let $G$ denote the group $SL_2(C)$ and let $\pi$ denote the projection map from $G$ onto the base manifold $\HI^3$ given by $\pi (g)=X=gg^*$. Then $\pi_*$  will denote the tangent map of $\pi$. It follows that $\pi_*(gA)=gAg^*$ for every left invariant vector field $g\rightarrow gA$ in $G$. The group $K$ acts on $G$ by right translations and the action is invariant and transitive on  each fiber $\pi^{-1}(g)$. Therefore $(G,\HI^3,\pi,K)$  is a principal K-bundle.

 . 
{\definition\label{horizdistr} The distribution 
spanned by the left-invariant vector fields $V(g)=gA$ with $A$ in $\fh$  and $g\in G$ shall be called horizontal and will be denoted by $\cal {H}$ . The integral curves of $\cal {H}$ in $G$ shall be called horizontal curves. Vertical curves are the integral curves of the left invariant distribution with values in $\fk$. A horizontal curve $g(s)$ is called a horizontal lift of a curve $X(s)$ in $\HI^3$ if $\pi(g(s))=X(s)$ for all $s$.}

{\theorem \label {subriem}(a). Every curve $X(s)$ in $\HI^3$ can be lifted to a horizontal curve $g(s)$. For any  horizontal lifts $g_1(s)$ and $g_2(s)$ of $X(s)$ there exists a vertical curve $a(s)$ such that $g_2(s)=g_1(s)a(s)$.

(b). If a horizontal lift $g(s)$ of a base curve $X(s)$ is given by  $\frac{dg}{ds}(s)=g(s)B(s)$ and $B(s)=\sum b_i(s)B_i$, then \[
||\frac{dX}{ds}(s)||^2 =4\langle B(s),B(s)\rangle =4( b_1^2(s)+b_2^2(s)+b_3^2(s)).\]}
where 
\[||\frac{dX}{ds}||^2=-\frac{dx_0}{ds}^2+\frac{dx_1}{ds}^2+\frac{dx_2}{ds}^2+\frac{dx_3}{ds}^2=-(\frac{dx}{ds},\frac{dx}{ds})_h\] 
\prf Let $g(s)$ denote any curve in $G$ that projects onto a given curve $X(s)$ in $\HI^3$  and let $A(s)$ and $B(s)$ denote curves in $\fg$ with $A(s)\in \fk$ and $B(s)\in \fh$ such that \[\frac{dg}{ds}(s)=g(s)(B(s)+A(s)).\]

Let $g_0(s)=g(s)\phi ^{-1}(s)$ where  $\phi (s)$ denote any solution of \[\frac{d\phi }{ds}(s)=\phi (s)A(s).\] Then, $g_0(s)$ projects onto $X(s)$ because $\phi (s)$ is a vertical curve, and furthemore, \[\frac{dg_0}{ds}(s)=g_0(s)(\phi (s)B(s)\phi ^{-1}(s)).\]
The relations (\ref{cartan}) imply that $\phi (s)B(s)\phi ^{-1}(s)$ belongs to $\fh$, and therefore $g_0 (s)$ is a horizontal curve. This proves the first statement in (a). The second statement in (a) follows directly from the definition of a principal bundle.

To prove part (b), let $g(s)$ ba a horizontal curve that projects onto $X(s)$ such that
\[\frac{dg}{ds}(s)=g(s)B(s)\]
Since $X(s)=g(s)g^*(s)$, $\frac {dX}{ds}(s)=2g(s)B(s)g^*(s)$.

Then Lemma \ref{hyptr} implies that
\bean(||\frac{dX}{ds}(s)||^2))I&=&-(\frac{dX}{ds}(s))(\frac{dX}{ds}(s))^{\dagger}=-4(g(s)B(s)g^*(s))(g(s)B(s)g^*(s)^{\dagger})\\&=&-4(g(s)B(s)g^*(s)(g^*(s)^{\dagger}B^{\dagger}(s)g^{\dagger}(s)=-4B(s)B^{\dagger}(s)\\&=&4(\langle B(s),B(s)\rangle)I
\eean

Therefore,
\[||\frac{dX}{ds}(s)||^2)=4\langle B(s),B(s)\rangle)\]

\qed
{\definition The hyperboloid $\HI^3$ shall be considered a Riemannian manifold with its metric given by
$\half\sqrt{\frac{dx_0}{ds}^2+\frac{dx_1}{ds}^2+\frac{dx_2}{ds}^2+\frac{dx_3}{ds}^2}$.}

Then previous lemma  implies
{\corollary\label{arclgth} Horizontal curves $g(s)$ defined by $\frac{dg}{ds}(s)=g(s)B(s)$ project onto curves $X(s)$ parametrized by arc-length, i.e. $||\frac{dX}{ds}(s)||=1$, if and only if $\langle B(s),B(s)\rangle =1$.} 
 
For the remainder of the paper all curves $g(s)$ in $G$ will be anchored at the identity ($g(0)=I$), and due to  the relation (\ref{hisom}) will be identified with  the curves in the orthonormal frame bundle of $\HI^3$. As in the spherical case

{\definition\label{hypdarb}  Curves $g(s)$ in $G$ will be called Darboux if $v_1(s)=g(s)B_1g^*(s)$ is equal to the tangent vector $\frac{dX}{ds}(s)$ of the projected curve $X(s)=\pi (g(s))$.}

The above condition means that Darboux curves are the solutions of \[\frac{dg}{ds}(s)=g(s)(B_1+A(s))\] 
for arbitrary curves $A(s)$ in $\fk$.
{\definition \label{horhypdarb}Horizontal curves $g(s)$ will be called horizontal anchored Darboux curves, or simply horizontal Darboux, if $g(0)=I$ and \[g^{-1}(s)\frac{dg}{ds}(s)=\Lambda (s)=R(s)B_1R^*(s)\]
 for some curve $R(s)$ in $K$ with $R(0)=I$.}

Since $K$ acts transitively on the sphere $S^3$, the horizontal Darboux curves could have been defined equivalently as the solutions $g(s)$ of the differential equation \[\frac{dg}{ds}(s)=g(s)\Lambda (s)\]
 such that $g(0)=I$, subject to further conditions that
 $\Lambda (s)\in \fh$, $\Lambda (0)=0$, and  $\langle \Lambda (s),\Lambda (s)\rangle =1$ for all $s\in [0,L]$.
{\definition The space of all hyperbolic Darboux curves and horizontal Darboux curves  shall be denoted by $\cD_h(L)$ and $Horiz(\cD_h)(L)$. }

{\definition\label{covderhyp} The covariant derivative  $\frac{D_g}{ds}(v)$ of a curve of tangent vectors $g(s)V(s)$, $V(s)\in \fb$, along a horizontal curve $g(s)$ in $G$, is defined by $\frac{D_g}{ds}(gV)(s)=g(s)\frac{dV}{ds}(s)$ for all $s\in [0,L]$.}

This notion of covariant derivative for vectors in the horizontal distribution $\cal{H}$ coincides with the usual notion of covariant derivative in the base manifold $\HI^3$ in the sense that\[\frac {D_{\pi(g)}}{ds}(\pi_*(gV))(s)=\pi_*((g(s)\frac{dV}{ds}))(s)\]

 The proof is simple and goes as follows:

 The covariant derivative in $\HI^3$ inherited from the ambient manifold $\IR^4$ is
  analogous to the spherical case, except that the Euclidean inner product in $\IR^4$  is replaced by the Lorentzian inner product, and is given by the following formula:
\[
\frac{D_x}{ds}(v)(s)=\frac{dv}{ds}(s)+{1\over 4}(\frac{dx}{ds}(s),v(s))_h x(s)
\]
{\bf{Remark}} The factor ${1\over 4}$ appears because the metric in this paper is equal to the half of the standard hyperbolic metric.
 
To adapt this formula to the representation by Hermitian matrices, first  note that if $V$ and $\Lambda $ are any Hermitian matrices in $\fg$ then
\[V^{\dagger}=-V\mbox{ and }\Lambda^{\dagger}=-\Lambda.\]

Then if $\pi_*(g(s)V(s)) $ is a curve of tangent vectors  defined along a curve $X(s)=\pi (g(s))=g(s)g^*(s)$ corresponding to a horizontal curve $g(s)$ in $G$ that is a solution of\[\frac{dg}{ds}(s)=g(s)\Lambda (s),\] 
  then the projected vectors $v$ and $\lambda$ are given by $v=\pi_*(gV)=2gVg^*$ and $\lambda=\pi_*(g\Lambda)=2g\Lambda g^*$. The results of Lemma \ref{hyptr} show that
\[
(v,\lambda)_hI=-\langle \Lambda , V\rangle I
\]

 It then follows that
\begin{eqnarray*}\frac{D_X}{ds}(\pi_*(g(s)V(s)))(s)&=&2\frac{d}{ds}(g(s)V(s)g^*)(s)-4\langle V,\Lambda\rangle
(g(s)g^*(s))\\&=& 2g(s)(\Lambda (s)V(s)+V(s)\Lambda (s))g^*(s)+2g(s)\frac{dV}{ds}(s)g^*(s)- (\langle V,\Lambda\rangle)g(s)g^*(s)\\&=&4\langle V,\Lambda\rangle)g(s)g^*(s)+2g(s)\frac{dV}{ds}(s)g^*(s)-4\langle V,\Lambda\rangle)g(s)g^*(s)\\&=&2g(s)\frac{dV}{ds}(s)g^*(s)
\end{eqnarray*}
\qed
 {\theorem\label{mixhypcovder} Suppose that $g(s,t)$ is a field of horizontal curves in $G$.  
 Let $V(s,t)=g(s,t)^{-1}\frac{\partial g}{\partial t})(s,t)$ and $\Lambda (s,t)=g(s,t)^{-1}\frac{\partial g}{\partial s})(s,t)$. Then, \[\frac{\partial \Lambda}{\partial t}(s,t)=\frac{\partial V}{\partial s}(s,t)\]
 and 
\[\frac{D_g}{ds}(\frac{\partial g}{\partial t})(s,t)=\frac{D_g}{dt}(\frac{\partial g}{\partial s})(s,t)\]}

\prf Let $\pi (g(s,t))=X(s,t)$. Then, as in any Riemannian manifold,
\[\frac{D_X}{ds}(\frac{\partial X}{\partial t})(s,t)=\frac{D_X}{dt}(\frac{\partial X}{\partial s})(s,t)\]
 The calculation  that preceeds the statement of the theorem shows that
\[\frac{D_X}{ds}(\frac{\partial X}{\partial t})(s,t)=2g(s,t)\frac{\partial V}{\partial s}(s,t)g^*(s,t)\mbox { and }\frac{D_X}{dt}(\frac{\partial X}{\partial s})(s,t)=2g(s,t)\frac{\partial \Lambda}{\partial t}(s,t)g^*(s,t)\]
Hence,
\[2g(s,t)\frac{\partial V}{\partial s}(s,t)g^*(s,t)=2g(s,t)\frac{\partial \Lambda}{\partial t}(s,t)g^*(s,t)\]
and consequently,\[\frac{\partial V}{\partial s}(s,t)=\frac{\partial \Lambda}{\partial t}(s,t)\]
The remaining statement follows directly from Definition \ref{covderhyp}
\qed
{\theorem \label{hypcurv}The curvature $\kappa (s)$ of the projection curve $X(s)=\pi (g(s))$ that corresponds to a horizontal curve $g(s)$ conforms to the following equation\[||\kappa(s)||^2=\langle \frac{d\Lambda}{ds}, \frac{d\Lambda}{ds}\rangle\] where $\Lambda (s)=g^{-1}(s)\frac{dg}{ds}(s)$.}

\prf \bean ||\kappa (s)||^2&=&||\frac{D_X}{ds}((\frac{dX}{ds})(s)||^2\\&=&||\pi_*(\frac{D_g}{ds}(\frac{dg}{ds}))||^2=||\pi_*(g(s) \frac{d\Lambda}{ds})||^2\\&=&\langle \frac{d\Lambda}{ds}, \frac{d\Lambda}{ds}\rangle\eean
\qed

When a horizontal curve $g(s)$ is a Darboux curve then its tangent vector $\Lambda(s)$ is  of the form\[\Lambda (s)=R(s)B_1R^*(s)\]
 for some curve $R(s)$ in $K$ that originates at the identity for $s=0$. If $R(s)$ is a solution of \[\frac{dR}{ds}(s)=R(s)A(s)=\sum u_i(s)A_i, then\]
\[\frac{d\Lambda}{ds}(s)=R(s)[B_1,A(s)]R^*(s)=R(s)(u_3(s)B_2-u_2(s)B_3)R^*(s)\], and therefore\[||\kappa (s)||^2=u_2^2(s)+u_3^2(s).\]

 The Serret-Frenet frame of the projected curve $X(s)=\pi (g(s))$, defined by the tangent vector $T(s)$, the normal vector $N(s)$ and the binormal vector $B(s)$ evolves according  to the well known equations:
\[\frac{D_X}{ds}(T)(s)=\kappa (s)N(s),\frac{D_X}{ds}(N)(s)=-\kappa (s)T(s)+\tau (s)B(s),\frac{D_X}{ds}(B)(s)=-\tau (s)N(s).\]
  
 When represented by  the equations (\ref{hypframe}) the Serret-Frenet vectors take the following form:\bean
 T(s)&=&\half R(s)E_1R^*(s)=R(s)B_1R^*(s)\\
N(s)&=&\half R(s)E_2R^*(s)=R(s)B_3R^*(s)\\
B(s)&=&\half R(s)E_3R^*(s)=-R(s)B_2R^*(s).
\eean

An easy calculation shows that\eq\label{curvtor}
u_1(s)=\tau (s),\,u_2(s)=-\kappa (s),\,u_3(s)=0.
\ee

\subsection{The symplectic form}

Similar to the spherical case, the space of anchored  hyperbolic horizontal Darboux curves will be considered  an infinite dimensional Fr\'echet manifold.

{\definition\label{perhypdarb} An anchored  horizontal Darboux curve $g(s)$ will be called frame-periodic if the tangent vector $\Lambda (s)$ is periodic with period $L$. The set of frame-periodic anchored horizontal Darboux curves will be denoted by $Horiz{(\cal PD}_h)(L).$}

Frame-periodicity requires periodicity of the curve $R(s)$ which in turn, imposes periodicity of its tangent vector $A(s)=\sum U_i(s)A_i$. The projections $X(s)$ in $\HI^3$ of frame periodic Darboux curves necessarily have periodic curvature and torsion, but need not be closed. 

On the other hand, smoothly closed curves $X(s)$ in $\HI^3$  generate periodic Serret-Frenet frames, and therefore can   be lifted always to periodic Darboux curves $g(s)$. If $g(s)$ is  the periodic solution of \[\frac{dg}{ds}(s)=g(s)(B_1+\sum u_i(s)A_i)\] with g(0)=I,
 then the corresponding horizontal curve $h(s)=g(s)R^*(s)$ is periodic if and only if $R(s)$ is a periodic solution of \[\frac {dR}{ds}(s)=R(s)(\sum u_i(s)A_i).\] 
For curves on the sphere and in a Euclidean space $R(s)$ can be identified with an orthonormal frame along the curve, which implies that smoothly periodic curves can be always lifted to periodic horizontal Darboux curves. For curves on the hyperboloid
 however, such a  statement requires a proof, since $R(s)$ does not  necessarily coincide with an orthonormal frame along the curve. 

{\theorem \label{hyptgt} (a) The tangent space ${T_g(Horiz(\cal D}_sh)(L)$  at an anchored  horizontal Darboux curve
$(g(s)$ consists of  tangent curves $v(s)=g(s)V(s)$ with $V(s)$ a Hermitian curve that  satisfies the following conditions: 
\eq\label{hypvariat} V(0)=0,\,\frac{dV}{ds}(0)=0,\,\langle \Lambda (s),\frac{dV}{ds}(s)\rangle=0\ee
where,\[\frac{dg}{ds}(s)=g(s)\Lambda (s).\] 

(b) For anchored frame-periodic horizontal Darboux curves  $v(s)=g(s)V(s)$  is a tangent vector at $g(s)$ if in addition to the properties in (a) the curve $\frac{dV}{ds}(s)$ is smoothly periodic having the period equal to $L$.}

\prf Let $h(s,t)$ denote a family of  anchored  horizontal Darboux curves such that $h(s,0)=g(s)$. Then $v(s)=\frac{\partial h}{\partial t}(s,t)_{t=0}$ is a tangent vector at $g(s)$. Since the curves $h(s,t)$ are anchored, $v(0)=0$.

Let $Z(s,t)$ and $W(s,t)$ denote the matrices defined by 
\[Z(s,t)=h(s,t)^{-1}\frac{\partial h}{\partial s}(s,t)\mbox {, }W(s,t)=h(s,t)^{-1}\frac{\partial h}{\partial t}(s,t).\]
It follows that $\Lambda (s)=Z(s,0)$ and $v(s)=g(s)V(s)$ with  $V(s)=W(s,0)$. Then $v(0)=0$ implies that $V(0)=0$. Furthemore,\[
\frac{\partial Z}{\partial t}(s,t)=\frac{\partial W}{\partial s}(s,t)\]

as a consequence of Lemma \ref{mixhypcovder}.   For $t=0$ the above equation reduces to
\[
\frac{dV}{ds}(s)=\frac{\partial W}{\partial s}(s,0)=U(s).
\]

 Since the curves $s\rightarrow Y(s,t)$ are Darboux for each $t$,\[\langle Z(s,t),Z(s,t)\rangle=1\mbox{, and } Z(0,t)=B_1\]
Therefore,
\[\langle \frac{\partial Z}{\partial t}(s,t),Z(s,t)\rangle=0\mbox{, and } \frac{\partial Z}{\partial t}(0,t)=0\]
which implies that $\langle \Lambda (s),U(s)\rangle=0$ and $U(0)=0$. 

Conversely any curve $V(s)$ in $\fh$ that satisfies (\ref{hypvariat}) can be realized by the perturbations $h(s,t)$ defined in the first part of the proof. The argument is similar to the one used for the spherical Darboux curves and goes as follows.

Let $V(s)$ be a curve specified by (\ref{hypvariat}), and let $\frac{dV}{ds}(s)=U(s)$.

  Define
\[
Z(s,t)={1\over {1+\phi ^2 (t)\langle U(s),U(s)\rangle }}(\Lambda (s)+\phi (t)U(s))
\]
 where $\phi (t)$ denote any smooth function such that $\phi (0)=0$ and $\frac{d\phi}{dt}(0)=1$.
Evidently, $Z(0,t)=B_1$ for all $t$, and an easy calculation shows that $\langle Z(s,t),Z(s,t)\rangle =1$. Therefore the solution of \[\frac{\partial h}{\partial s}(s,t)=h(s,t)Z(s,t)\] with $h(0,t)=I$ corresponds  to an anchored  horizontal Darboux curve for each $t$. Since $U(s)=\frac{\partial Z}{\partial t}(s,0)$ the proof of part (a) is finished.

To prove part (b) assume that the curves $h(s,t)$ used in part (a) belong to $Horiz({\cal {PD}}_s )(L)$. Then curves $s\rightarrow Z(s,t)$ are periodic with period $L$ for each $t$. Since curves $s\rightarrow \frac{\partial Z}{\partial t}(s,t)(s,t)$ preserve this periodicity, $U(s)=\frac{\partial Z}{\partial t}(s,0)$ is periodic with period $L$.

\qed
 
{\definition \label{hypsymp} The symplectic form $\omega $ on the  space of anchored horizontal hyperbolic Darboux curves  is defined as follows\eq\omega_{g\Lambda}(V_1,V_2)={1\over i}\int_0^L\langle \Lambda (s),[\frac{DV_1}{ds}(s),\frac{dV_2}{ds}(s)]\rangle \,ds\ee for any tangent vectors $V_2(s),V_2(s)$ at a horizontal curve $g(s)$.}

{\bf{Remark}}  As it stands the above definition may be problematic, since the issues of non-degeneracy and closedness of $\omega$ have not been dealt with yet. However, such questions are easily removed because the hyperbolic form $\omega$ is isomorphic to its spherical analogue. The proof is as follows:

Let \[\tilde{\Lambda}=i\Lambda ,\,\tilde{U_1}=i\frac{dV_1}{ds},\,\tilde{U_2}=i\frac{dV_2}{ds}.\]
It follows from (\ref{cartan}) and  Table 1 that $\tilde{\Lambda},\,\tilde{U_1},\,\tilde{U_2}$ all belong to $\fk$ and that$\langle \tilde{\Lambda},\tilde{U_i}\rangle=0$ for $i=1,2$.
Therefore, 
\bean \omega_{\Lambda}(V_1,V_2)&=&{1\over i}\int_0^L\langle \Lambda (s),[\frac{DV_1}{ds}(s),\frac{dV_2}{ds}(s)]\rangle \,ds\\
&=&-\int_0^L\langle \tilde{\Lambda }(s),[\tilde{U_1}(s),\tilde{U_2}(s)]\rangle \,ds=\tilde{\omega}_{\tilde{\Lambda}}(\tilde{V_1},\tilde{V_2})\eean
under the identification of the curve $X(s)$ in $K$ with the solution of\[\frac{dX}{ds}(s)=X(s)\tilde{\Lambda}(s)\] that satisfies $X(0)=I.$ 

It follows that $\tilde{\omega}$ coincides with the symplectic form for the spherical horizontal Darboux curves. The isomorphism, apart from justifying the choice of sign in (\ref{dsymp}) also makes transparent the proof of  
{\theorem  Both $Horiz({\cal D}_h)(L)$ and $Horiz({\cal PD}_h)(L)$ are symplectic manifolds relative to $\omega$ defined by (\ref{hypsymp}).}
{\remark 
It may be appropriate to point out that both the spherical and the hyperolic 
 symplectic form  in this paper may be seen as natural adaptations to curves of the  standard symplectic form $\omega$ 
 on the sphere $S^2$  given explicitly by 
\[ \omega_{\gamma }(a,b)={\gamma }\cdot {(a\times b})\] where $a$
an $b$ are tangent vectors at a point $\gamma $ on $S^2$. In the light of this observation,
the symplectic structure of anchored curves with periodic frames is isomorphic 
to the symplectic structure of anchored loops on the sphere.}

In each of the two cases 
the right action of $SU_2$ extends to the space of horizontal anchored Darboux curves with $(g,a)\rightarrow ag(s)a^\ast$ for each
 horizontal curve $g(s)$ and each $a\in SU_2$. This action is symplectic relative to the forms  used in this paper
and $J(g)=\int_0^L\Lambda(s)\,ds$ is the moment map associated 
with this action (under the implicit assumption that the dual $\fg^\ast$ of the Lie algebra $\fg$ is identified with $\fg$ 
via the trace form ).
 
 The moment map induces a function 
$J_A(g)=\int_0^L\leangle \Lambda (s),A \riangle\,ds$
on the space of  horizontal anchored curves for each  element $A\in\fg$. The Hamiltonian vector field induced by this function coincides 
with the 
infinitesimal generator of the action-induced one-parameter group of transformations 
$\{e^{tA}\gamma e^{-tA}\}$.
Then it is well known ([1]) that $J$ is an integral of motion for each Hamiltonian function 
which is
 invariant under the action. 

The moment map will be taken up again in the problems of mathematical physics further down in the text.
There is another symplectic form on the space of anchored curves given by the following expression:
\[
\Omega_{g\Lambda}(V_1,V_2)=\int_0^L \leangle \Lambda ,[V_1,V_2]\riangle\, ds 
\]

Such a form is mentioned elsewhere in the literature (see for instance [2],[3], and [13]).
This form is compatible in the sense of Magri ([15]) with the form used in this paper
 and can  be used to get the 
integrability results
for systems which are bi-Hamiltonian (as outlined in ([2]). However, such investigations seem too particular for the scope of this paper and will not be pursued here.

\section{The Hamiltonian flow of $\half\int_0^L{k^2(s})\,ds$}

  Horizontal Darboux curves $g(s)$ are parametrized by the matrices $\Lambda (s)=R(s)E_1R^*(s)$ with curves $R(s)$  the solutions of \[\frac{dR}{ds}(s)=R(s)\sum u_i(s)A_i\]
in $K$ with $R(0)=I$. The matrix $E_1$ is equal to $A_1$ in the spherical case and equal to $B_1$ in the hyperbolic case.

Theorems \ref{curvdarb} and \ref{hypcurv} show that $||\frac{d\Lambda }{ds}(s)||^2=\kappa^2(s)$ where $\kappa (s)$ is the curvature of the projected curve in the base manifold. Hence the function
\[f(g)=\half\int\limits_0^L||\frac{d\Lambda }{ds}(s)||^2\,ds\]
on the space of anchored Darboux curves may be also seen as the function
  \[ f(x)=\half\int\limits_0^L\kappa ^2(s)\,ds\]
on the space of  projected curves $x(s)=\pi (g(s))$. 

As a function on the space of horizontal Darboux curves $f$ defines a 
 Hamiltonian vector field  $\mathcal{X}_f$  through the symplectic formalism  described in the first part of the paper. Remarkably, $ \mathcal{X}_f$ leads to  Heisenberg's magnetic equation and the non-linear Schroedinger's equation when $f$ is restricted to the frame-periodic Darboux curves. A derivation of this fact, together with  the connections to the known results in the literature constitute the subject matter for the remaining part of the paper.

\subsection{ Heisenberg's magnetic equation} 

Although conceptually alike, the calculations  in the spherical setting are different in several aspects  from those in the hyperbolic setting and will be done separately in each of the above mentioned cases.
 
{\bf{Hyperbolic Darboux curves.}} To calculate the directional derivative $df_{\Lambda}(V)$, let $\hat{g}(s,t)$ be a family of anchored horizontal Darboux curves  that are the solutions of \[\frac{\partial\hat{g}}{\partial s}=\hat{g}(s,t)\hat{\Lambda}(s,t)\] such that \[\hat{g}(s,0)=g(s),\,\hat{\Lambda}(s,0)=\Lambda (s),\,\frac {\partial \hat {\Lambda}}{\partial
 t}(0,s)=\frac{dV}{ds}(s)\]. 

The directional derivative $df_{\Lambda}(V)$   is given by \[
df_{\Lambda}(V)=\half{\partial\over\partial t}\int_0^L\leangle{\partial\hat{\Lambda}\over\partial s}
		(s,t),{\partial\hat{\Lambda}\over\partial s}
		(s,t)\riangle\,ds|_{t =0.}
\]
It follows that

\begin{eqnarray*}
\half{\partial\over\partial t}\int_0^L\leangle{\partial\hat{\Lambda}\over\partial s}
		(s,t),{\partial\hat{\Lambda}\over\partial s}
		(s,t)\riangle\,ds|_{t =0}
	& = &\int\limits_0^L\leangle{\partial\hat\Lambda\over
		\partial s}(s,t),\frac{\partial}{\partial s}
		{\partial\hat\Lambda\over\partial t}(s,t)
		\riangle\,ds|_{t=0}\\
	& = &\int\limits_0^L\leangle
		{d\Lambda\over ds},\frac{d}{ds}
		\left({dV\over ds}\right)\riangle\, ds\\
	& = &-\int\limits_0^L
		\leangle {d^2\Lambda\over ds^2},
		{dV\over ds}\riangle\, ds+	
		\leangle{d\Lambda\over ds},{dV\over ds}
		\riangle |_{s=0}^{s=L.}
\end{eqnarray*}

  In the space of  frame-periodic horizontal Darboux curves  the boundary terms $\leangle{d\Lambda\over ds},{dV\over ds}
		\riangle |_{s=0}^{s=L}$ are equal to $0$ because of periodicity.
  Consequently, 
\[ df_\Lambda (V)= -\int\limits_0^L\leangle {d^2\Lambda\over ds^2},
{dV\over ds}\riangle\, ds .\]
The Hamiltonian vector field  is of the form $\mathcal {X}_f(g)=gF$  for some Hermitian matrix $F(s)$ that satisfies\[ df_\Lambda (V)= 
\frac{1}{i}\int_0^L\leangle \Lambda(s),
[{dF\over ds},{dV\over ds}]\riangle\, ds\] for an arbitrary tangential direction $V(s)$.
The above is equivalent to\eq\label{int}\int\limits_0^L\leangle ({d^2\Lambda\over ds^2}+\frac{1}{i}[\Lambda(s),
{dF\over ds}]),{dV\over ds}\riangle\, ds=0.\ee

It then follows from Theorem \ref{hyptgt} that $V(s)$ can be taken as $V(s)=i\int_0^s[\Lambda (t),U(s)]\,dt$ where $U(s)$ is an arbitrary curve of Hermitian matrices with $U(0)=0$.
 Then equation (\ref{int}) becomes\[ 
\int\limits_0^L\leangle (i[{d^2\Lambda\over ds^2},\Lambda (s)]+[[\Lambda(s),
{dF\over ds}],\Lambda (s)]),U(s)\riangle\, ds=0.\]
It follows that \eq\label{vectfld}i[{d^2\Lambda\over ds^2},\Lambda (s)]+[[\Lambda (s),
{dF\over ds}],\Lambda (s)]=0\ee because $U(s)$ is sufficiently arbitrary. 
 The following Lie bracket identity is true\[[[A,B],C]=\langle B,C\rangle A-\langle A,C\rangle B\]  for any Hermitian matrices $A,B,C$ as can be readily verified through Table 1. 

Therefore,\[[\Lambda (s),[\Lambda(s),
{dF\over ds}]]=-{dF\over ds}\]
and equation (\ref{vectfld}) becomes \[i[{d^2\Lambda\over ds^2},\Lambda (s)]-{dF\over ds}=0\]
It follows that \[\mathcal{X}_f(g)=g(s)F(s)\mbox { with }F(s)=i\int_0^s [{d^2\Lambda\over dx^2}(x),\Lambda (x)]\,dx\]
 is the Hamiltonian vector field that corresponds to  $f$. 

The integral curves $t\rightarrow g(s,t)$  of  $\mathcal{X}_f$ are the solutions of the following partial differential equations \begin{eqnarray}\label{hamflow}
{\partial g\over \partial t}(s,t)& = &g(t,s)i\int_0^s [{d^2\Lambda\over dx^2}(x,t),\Lambda (x,t)]\,dx\\
{\partial g\over\partial s}(s,t)& = &g(t,s)\Lambda(t,s).
\end{eqnarray}

 Theorem \ref{mixhypcovder} implies that the matrices $\Lambda (s,t)$ evolve according to 
\begin{equation}\label{heisen}
{\partial\Lambda\over\partial t}(t,s)=
i\left[{\partial^2\Lambda\over\partial s^2} ,\Lambda (t,s)\right].
\end{equation}

Equation (\ref{heisen} when expressed in terms of the coordinates  $\lambda(s,t)$ of $\Lambda (s,t)$ 
relative to the basis of the Hermitian Pauli matrices becomes: 

\begin{equation}\label{heiseneuc}
\frac{\partial \lambda}{\partial t}(t,s)=\lambda (t,s)\times \frac{\partial^2 \lambda }
{\partial s^2}(s,t).
\end{equation}

{\remark Equation (\ref{heiseneuc}) is well known in the literature in applied mathematics.
  L.D.Faddeev and L.A. Takhtajan refer to it as 
the continuous isotropic
 Heisenberg feromagnetic model([5], Part II.,Chapter 1) which they treat  in an ad hoc manner as an equation in the space of Hermitian matrices. V.I. Arnold and B.Khesin ([2]) connect (\ref{heiseneuc}) to the filament equation  which they furher consider as a special type of a Landau-Lifschitz equation on $so_3(R)$.}

{\bf{Spherical Darboux curves.}} 
 The derivation of the corresponding Hamiltonian equations on the sphere is similar to 
the preceding
 case except for the details related to the covariant derivative. Recall that the tangent space ${T_X(Horiz(\cal D}_s )(L)$  at an anchored  horizontal Darboux curve
$X(s)$ consists of  tangent curves $v(s)=X(s)V(s)$ with $V(s)$ the solution of
\begin{equation}\label{displ}
\frac{dV}{ds}(s)=[\Lambda (s),V(s)]+U(s)
\end{equation}
 with $V(0)=0$. The matrix $\Lambda (s)$ is the tangent vector of $X$, i.e.,
 \[\frac{dX}{ds}(s)=X(s)\Lambda (s)\]
 and $U(s)$ is a
curve in $\fk$ subject to $U(0)=0$ and $\langle \Lambda (s),U(s)\rangle=0$.

Let $v(s)=X(s)V(s)$ be a fixed tangent vector at an anchored horizontal Darboux curve $X(s)$. To find the appropriate expression for $df_X(V)$ the directional derivative of $f$ at $X$ in the direction $V$, let $Y(s,t)$ denote a family of  anchored  horizontal Darboux curves such that $Y(s,0)=X(s)$ and such that $v(s)=\frac{\partial Y}{\partial t}(s,t)_{t=0}$. 

Let $Z(s,t)$ denote the matrices defined by 
\[\frac{\partial Y}{\partial s}(s,t)=Y(s,t)Z(s,t)\]
It follows that $\Lambda (s)=Z(s,0)$, and that $V(s)$ is the solution of (\ref{displ}) with $U(s)=\frac{\partial Z}{\partial t}(s,0)$.
 
Then,\bean df_\Lambda(V)&=&\half\frac{\partial}{\partial t}\int_0^L\langle \frac{\partial Z}{\partial s}(s,t),\frac{\partial Z}{\partial s}(s,t)\rangle \,ds|_{t=0}\\
&=&\int_0^L\langle \frac{d\Lambda}{ds} (s),\frac{dU}{ds}(s)\rangle \,ds\\
&=&-\int_0^L\langle \frac{d^2\Lambda}{ds^2}(s),U(s)\rangle\,ds+\langle \frac{d\Lambda}{ds}(s),U(s)\rangle_{s=0}^{s=L}.\eean

 Analogous to the hyperbolic case the boundary terms vanish in the frame-periodic case, and therefore 
\eq\label{1} df_X(V)=-\int_0^L \leangle \frac{d^2 \Lambda}{ds^2},U(s)\riangle \,ds
\ee
The Hamiltonian vector field $\mathcal{X}_f$ that corresponds to $f$ is  of the form \[\mathcal{X}_f(X)(s)=X(s)F(s)\] for some curve $F(s)\in \fk$. Since $\mathcal{X}_f(X)\in T_X Horiz(\mathcal{PD})_s(L)$ $F(s)$ is the solution  of \[\frac{dF}{ds}(s)=[\Lambda (s),F(s)]+U_f(s)\]with $F(0)=0$  for some curve $U_f(s)\in\fk$ that satisfies\[U_f(0)=0\,\mbox{ and }\langle \Lambda (s),U_f(s) \rangle =0.\]

The curve $U_f(s)$ is determined by the symplectic form $\omega$ in Definition \ref{symp} through the usual relation\eq \label{2}
df_\Lambda (U)=-\int_0^L\langle \Lambda (s),[U_f(s),U(s)]\rangle\,ds .\ee
The curve $U(s)$ satisfies $U(0)=0$ and $\langle \Lambda (s),U(s)\rangle =0$, and is otherwise arbitrary.
 
Equations (\ref{1}) and (\ref{2}) yield
\eq\label{3}
\int_0^L\leangle \frac{d^2 \Lambda}{ds^2}-[\Lambda (s),U_f(s)],U(s)\riangle \,ds=0.
\ee
 The curve $U(s)$ can be written as $U(s)=[\Lambda (s),C(s)]$ where $C(s)$ is any curve that satisfies $C(0)=0$. In that case equation (\ref{3}) becomes
\eq\label{4}
\int_0^L\leangle [\frac{d^2 \Lambda}{ds^2},\Lambda ] -[[\Lambda (s),U_f(s)],\Lambda ],C(s)\riangle \,ds=0.
\ee
Since $C(s)$ is arbitrary,
\eq\label{5}
[\frac{d^2 \Lambda}{ds^2},\Lambda ] -[[\Lambda (s),U_f(s)],\Lambda ]=0.\ee
Lemma \ref{triplebracket} implies that\[[[\Lambda (s),U_f(s)],\Lambda ]=U_f\]
and therefore,\[U_f=-[\Lambda ,\frac{d^2 \Lambda}{ds^2}].\]
 The  integral curves $t\rightarrow X(s,t)$ of  the vector field $\mathcal{X}_f$  are the solutions of\begin{equation}
\frac{\partial X}{\partial t}(s,t)=X(s,t)F(s,t)\mbox{, and } 
\frac{\partial X}{\partial s}
=X(s,t)\Lambda(s,t)
\end{equation}
where $F(s,t)$ is the solution of \[\frac{\partial F}{\partial s}(s,t)=[\Lambda (s,t),F(s,t)]-[\Lambda (s,t),\frac{d^2 \Lambda}{ds^2}(s,t)].\]
But then according to Lemma \ref{mixcovder}
\[
\frac{\partial \Lambda}{\partial t}(s,t)-\frac{\partial F}{\partial s}(s,t)+
[\Lambda(s,t),F(s,t)]=0
\]
hence,
\begin{equation}\label{6}
\frac{\partial \Lambda}{\partial t}(s,t) =
-[\Lambda (s,t),\frac{\partial ^2 \Lambda}{\partial s^2}(s,t)]=[\frac{\partial ^2 \Lambda}{\partial s^2}(s,t),\Lambda (s,t)].
\end{equation}

  Equation (\ref{6}) describes the flow of the Hamiltonian vector field in the spherical case. 

The reader should keep in mind however, that in the hyperbolic case $\Lambda$ is Hermitian,
 while in the spherical case, $\Lambda$ is skew-Hermitian: to pass from the hyperbolic case to the spherical case
  multiply $\Lambda (t,s)$ in equation (\ref{heisen}) by $i$.
{\definition Equations (\ref{heisen}) and (\ref{6}) shall be called Heisenberg's magnetic equations.}  
\subsection{The non-linear Schroedinger equation}
  
   Each solution $\Lambda (s,t)$ of Heisenberg's magnetic equation is generated by a family of  periodic frames $t\rightarrow R(s,t)$ in $\fk$ through $\Lambda (s,t)=R(s,t)B_1R^*(s,t)$ in the hyperbolic case, and through $\Lambda (s,t)=R(s,t)A_1R^*(s,t)$ in the spherical case.
 Curves $R(t,s)$ then evolve according
 to the differential equations:
\[
\frac{\partial R}{\partial s}(s,t)=R(s,t)U(s,t)\mbox{, and  }\frac{\partial R}{\partial t}=R(s,t)V(s,t) 
\]
 for some curves of matrices $U(s,t)$ and $V(s,t)$ in $\fk$. Matrices $U(s,t)$ and $V(s,t)$ are not independent of each other since they conform to 
\begin{equation}\label{zerocurv}
	{\partial U\over\partial t}(s,t)-
	{\partial V\over\partial t}(s,t)+	
	[U(s,t),V(s,t)]=0
\end{equation}
  as demonstrated earlier in the paper ( Lemma {\ref{mixcovder}).
Moreover, $V(0,t)=0$ because the  horizontal Darboux curves are anchored at $s=0$. Equation (\ref{zerocurv}) then implies that $\frac{dU}{dt}(0,t)=0$.
{\theorem\label{h-s}  Let $U(s,t)=\sum u_j(s,t)A_j$ generate a solution of Heisenberg's magnetic equation and let $u(s,t)=u_2(s,t)+iu_3(s,t)$. Then, \[
\psi(s,t)=u(s,t)\exp{(i\int_0^s u_1(x,t)\,dx )}
\] 
 is a solution of the non-linear Schroedinger's equation
\begin{equation}\label{nls}
	{\partial\over\partial t}\psi (s,t)
	=i{\partial^2\psi\over\partial s^2}(s,t)
	+i\half|\psi (s,t)|^2\psi (s,t).
\end{equation}}

The following lemma, whose proof can be easily obtained from the Lie brackets in Table 1, is useful for the proof of the theorem.
{\lemma \label{triplelie}
\[[A,[A,B]]=\leangle A,B \riangle A-\leangle A,A \riangle B \]
for any $A,B \in \fh$.
\[ [[A,B],B]=\leangle B,B \riangle A-\leangle A,B \riangle B\]
for any $A\in su_2$ and $B\in\fh$.}

{\bf{Proof of the theorem}} The proof of the theorem will be done for the hyperbolic case although the arguments are the same in both cases, as will become clear below.

Since $\Lambda(s,t)=R(s,t)B_1R^*(s,t)$, it
follows that
\begin{eqnarray*} 
	{\partial\Lambda\over\partial t}& = &
		{\partial R\over\partial t} B_1 R^{\ast}+R B_1
		{\partial R^{\ast}\over\partial t}\\
	& = &R(VB_1-B_1V)R^{\ast}\\
	& = &R[B_1,V]R^{\ast}
\end{eqnarray*}

Similarly,
\[ {\partial\Lambda\over\partial s}=R[B_1,U]R^{\ast},
{\rm\quad and\quad}{\partial^2\Lambda\over\partial s^2}
=R\left(\left[[B_1,U],U\right]+\left[B_1,{\partial U\over\partial s}
 \right]\right)R^{\ast}\]

The fact that $\Lambda (s,t)$ evolves according to Heisenberg's magnetic equation implies that
\begin{equation}\label{heisenlie}
[B_1,V]=i([[[B_1,U],U],B_1]+[[B_1,\frac{\partial U}{\partial s}],B_1])
\end{equation}

According to the relations in Lemma \ref{triplelie}
\[ 
[[B_1,U],U]=\leangle U,B_1 \riangle U-\leangle U,U \riangle B_1 = -(u_2^2+u_3^2)B_1 +{u_1u_2}B_2+ {B_3u_1u_3}\]
, and therefore 
\[
[[[B_1,U],U],B_1]={u_1}{u_3}A_2-{u_1}{u_2}A_3.
\]

Similarly,
\[
[B_1,\frac{\partial U}{\partial s}]=\frac{\partial u_3}{\partial s}B_2-\frac{\partial u_2}{\partial s}B_3\mbox
{, and }[[B_1,\frac{\partial U}{\partial s}],B_1]= -\frac{\partial u_3}{\partial s}A_3 -
\frac{\partial u_2}{\partial s}A_2.
\]
Then equation (\ref{heisenlie}) reduces to \[ [B_1, V]=i(u_1(u_3A_2-u_2A_3)-\frac{\partial u_3}{\partial s}A_3-
\frac{\partial u_2}{\partial s}A_2)\] which can also be written as 
\[
[B_1,V]=-u_1({u_3}B_2-{u_2}B_3)+\frac{\partial u_3}{\partial s}B_3+\frac{\partial u_2}{\partial s}B_2
\]
because $A_j=iB_j, j=1,2,3$.

 Let $V(s,t)=v_1(s,t)A_1+v_2(s,t)A_2+v_3(s,t)A_3$.
 Then, $[B_1,V]={v_3}B_2-{v_2}B_3$ , which together with the 
relations above yields
\begin{equation}\label{v1}
v_2=-u_1u_2-\frac{\partial u_3}{\partial s}\mbox{, and } v_3=-u_1u_3 +\frac{\partial u_2}{\partial s}.
\end{equation}
Equation (\ref{v1}) written in terms of  the complex functions $u=u_2+iu_3$ and $v=v_2+iv_3$ becomes
\begin{equation}\label{v2}
v(s,t)=-u_1(s,t)u(s,t)+i\frac{\partial u}{\partial s}(s,t).
\end{equation}

The zero curvature equation ${\partial\U\over\partial t}-
{\partial V\over\partial s}+[U,V]=0$ implies that
\begin{equation}\label{u_1v_1}
	{\partial u_1\over\partial t}=
	{\partial v_1\over\partial s}+\half
	{\partial \over\partial s}(u_2^2+u_3^2)
\end{equation}
and that 
\begin{equation}\label{uv}
	{\partial u\over\partial t}= i{\partial^2u\over\partial s^2}-
			2u_1{\partial u\over\partial s}
			-{\partial u_1\over\partial s}u
			-i(v_1+u_1^2)u .
\end{equation}

Equation (\ref{u_1v_1}) implies that
\[ {\partial\over\partial t}\int\limits_0^s u_1(x,t) dx
=v_1(s,t)+\half(u_2^2(s,t)+u_3^2(s,t))+c(t)\]
for some function $c(t)$. However this function must be equal to zero because the Darboux  curves are anchored, and therefore
 $v_1(t,0)=0$ and $u(t,0)=0$.

Upon the substitution $v_1(s,t)={\partial\over\partial t}
\int\limits_0^su_1(x,t)dx-\half|u(s,t)|^2$  the equation (\ref{uv})
becomes
\begin{equation}\label{uv1}
{\partial u\over\partial t}+iu{\partial\over\partial t}
\int\limits^su_1(t,x)\,dx=i{\partial^2u\over\partial s^2}
-2u_1{\partial u\over\partial s}-u{\partial u_1\over\partial s}
-i\left(-\half|u|^2+u_1^2\right)u.
\end{equation}

After the multiplication by
$\exp{(i\int_0^s u_1(x,t)\,dx})$ equation (\ref{uv1}) can be expressed as
 \[
 {\partial \over\partial t}(u\exp{(i\int_0^su_1\,dx)})
=(i{\partial^2u\over\partial s^2}-2u_1{\partial u\over\partial s}
-u{\partial u_1\over\partial s}-i(u_1^2-\half|u|^2))
u\exp{(i\int_0^su_1\,dx)}\]
because \[{\partial\over\partial t}(u(s,t)\exp{(i\int_0^s u_1(x,t)\,dx)})=\exp{(i\int_0^s u_1(x,t)\,dx)}({\partial u\over\partial t}+iu{\partial\over\partial t}
{\int\limits_0^su_1(x,t)\,dx)}.\]

 The function \[\psi(s,t)=u(s,t)\exp{(i\int_0^s u_1(t,x)\,dx)}\]
 satisfies \[{\partial\psi\over\partial s}=\left(
{\partial u\over\partial s}+iuu_1\right)\exp{(i\int_0^s u_1(x,t)\,dx})\]
and
\[{\partial^2\psi\over\partial s^2}=
\left({\partial^2u\over\partial s^2}+2iu_1
{\partial u\over\partial s}+iu{\partial u_1\over\partial s}
-u_1^2u\right)\exp{i\int_0^s u_1(x,t)\,dx}.\]

It follows that
\[ i{\partial^2\psi\over\partial s^2}=
\left(i{\partial^2u\over\partial s^2}-2u_1
{\partial u\over\partial s}-u{\partial u_1\over
\partial s}-iu_1^2u\right)\exp{i\int_0^s u_1(x,t\,dx}\]
, and therefore
  \[
	{\partial\over\partial t}\psi (t,s)
	=i{\partial^2\psi\over\partial s^2}
	+i\half|\psi|^2\psi .
\]

\qed

In the spherical case the evolution along Heisenberg's magnetic equation leads to
\begin{equation*}
[A_1,V]=([[[A_1,U],U],A_1]+[[A_1,\frac{\partial U}{\partial s}],A_1]).
\end{equation*}
 The preceeding equation is the same as  equation (\ref{heisenlie}) because $A_1=iB_1$. Therefore the calculations that led to the non-linear Schroedinger equation in the hyperbolic case are equally valid in the spherical case with the same end result.\\

  The steps  taken in the passage from Heisenberg's equation to the Schroedinger's equation are reversible. Any solution $\psi(s,t)$ of (\ref{nls}) generates a solution of the zero-curvature equation ( as demonstrated in [5]). Simply  let \[U=\half
\begin{pmx}{cc}
	0 & \psi\\
	&\\
	-\overline{\psi} & 0
\end{pmx}\mbox{ and }V=\half
\begin{pmx}{cc}
	i{\partial\over\partial s}|\psi|^2 & i\psi(|\psi|^2+c)\\
	&\\
	i\overline{\psi}|\psi|^2 & -i{\partial\over\partial s}
		|\psi|^2
\end{pmx}.\] 
 Then curves $R(s,t)$ that evolve according
 to the differential equations:
\[
\frac{\partial R}{\partial s}(s,t)=R(s,t)U(s,t),\,\frac{\partial R}{\partial t}(s,t)=R(s,t)V(s,t) 
\]
generate the solutions  $\Lambda (s,t)$ of Heisenberg's magnetic equation through the familiar formulas $\Lambda (s,t)=R(s,t)B_1R^*(s,t)$ or $\Lambda (s,t)=R(s,t)A_1R^*(s,t)$ depending on the case.

To correlate the findings of this paper with the related existing literature, which almost exclusively deals with curves in $\IR^3$, it seems appropriate to include  a discussion of the only remaining  simply connected three dimensional symmetric space, namely the Euclidean space.

\subsection{ Euclidean Darboux curves }

 The semidirect product of $\fh$ with $K$ is the most convenient  setting for comparisons with non-Euclidean Darboux curves. Recall that the semidirect product $V\triangleright K$ of a vector space $V$  and a group $K$ which acts linearly  on $V$ consists of pairs $(x,R)$ with $x\in V$ and $R\in K$. The group operation is given by $(x,R)(y,T)=(x+Ry,RT)$ for any elements $(x,R)$ and $(y,T)$. 
The Lie algebra $V\triangleright \fk$ of  the semidirect product $V\triangleright K$ consists of pairs $(a,A)$ with $a\in V$ and $A\in\fk$ with the Lie bracket $[(a,A),(b,B)]=(A(b)-B(a),[A,B])$. Both the vector space $V$ and the Lie algebra $\fk$ can be embedded in the Lie algebra of the semidirect product via the embeddings $a\rightarrow (a,0)$ and $A\rightarrow (0,A)$. With this identification, $V\triangleright \fk=V\oplus\fk$ and\[
[V,V]=0,\,[V,\fk]=V,\,[\fk,\fk]=\fk.\]

The group $K=SU_2$ acts linearly on the space of Hermitian matrices $\fh$ by $R(x)=RxR^*$ for $x\in V$ and $R\in K$, and $A(a)=[a,A]$ for $a\in\fk$ and $A\in\fk$. Therefore, the Lie bracket in $\fh\triangleright \fk$ is given by  \[[(a,A),(b,B)]=([b,A]-[a,B],[A,B]).\] The space of Hermitian matrices endowed with the metric induced by the trace form becomes a three dimensional Euclidean space, that will be denoted by $E^3$, while the semidirect product of $K$ with $E^3$ will be denoted by $S_K (E_3)$ .
 
The group $S_K(E_3)$ acts on $E^3$ by $(x,R)(y)=R(y)+x$ for each $(x,R)\in S_K(E_3)$ and each $y\in E^3$.
The action is transitive, and $K$ is equal to the isotropy group of the orbit through the origin $  y=0$. The Euclidean space $E^3$, when identified with the orbit through the origin becomes   the coset space $S_K(E_3)/K$.
 The preceeding action extends to an action on the tangent bundle in which
 a tangent vector $v$ at $y$ is taken to the the tangent vector $R(v)$ at $x$ under the action by an element $(x,R)\in S_K(E_3)$. The action  on the tangent bundle extends further
to  an action on the orthonormal frame bundle of $E^3$ such that a frame $(v_1,v_2,v_3)$ at a point $y\in E^3$
 is taken to the frame $(R(v_1),R(v_2),R(v_3))$ at $x$ under the action by an element $(x,R)\in S_K(E_3)$. The kernel of this action consists of $\pm I$, and hence
  $S_K(E_3)/\{\pm I\}$ can be identified with the positively oriented orthonormal frame bundle of $E^3$ as the orbit through the standard frame $(B_1,B_2,B_3)$ at the origin. 

 In the left-invariant representation of the tangent bundle of $S_K(E_3)$, the tangent vectors at a point $(x,R)$  are given by pairs $(R(a),RA)$ with $a\in E^3$ and $A\in\fk$. Hence curves $(x(s),R(s))$ in $S_K(E_3)$ are represented by differential equations\eq\label{eucdiff}
\frac{dx}{ds}(s)=R(s)(a(s)),\,\frac{dR}{ds}(s)=R(s)A(s).
\ee
The terminology concerning  Darboux curves in non-Euclidean cases extends naturally to the Euclidean setting. In particular,
  curves $(x,R)\in S_K(E_3)$ are Euclidean Darboux curves if $\frac{dx}{ds}(s)=v_1(s)=R(s)(B_1)$, which holds whenever $a(s)=B_1$. Horizontal  Euclidean Darboux curves are the projections $x(s)$ of Euclidean Darboux curves. Anchored horizontal Darboux curves are the solutions of \[ \frac{dx}{ds}(s)=R(s)(B_1),\,x(0)=0\] with $R(s)$ an arbitrary curve in $K$ that originates at $I$ when $s=0$. Frame-periodic  Darboux curves $(x,R)$ conform to the periodicity of $R(s)$ with its period equal to the length of $x(s)$.

For any horizontal Darboux curve $x(s)$, \[\frac{d^2x}{ds^2}(s)=R(s)([B_1,A(s)])\]
, and therefore\[\kappa^2(s)=||\frac{d^2x}{ds^2}(s)||^2=u_2^2(s)+u_3^2(s)\]
where $A(s)=\sum u_i(s)A_i$.
 The frame $R(s)$  is a Serret-Frenet frame if $A(s)=\tau (s)A_1+\kappa (s)A_3$, in which case
the frame vectors  $T(s),N(s),B(s)$ are given by  \[T(s)=R(s)(B_1),\,N(s)=R(s)(B_2),\,B(s)=R(s)(B_3).\]

The reader may easily verify that the tangent space at each anchored horizontal Darboux curve $x(s)$ consists of curves $v(s)$ such that 

(a) $v(0)=\frac{dv}{ds}(0)=0$, and

(b)  $ \langle \frac{dx}{ds}(s),\frac{dv}{ds}(s)\rangle=0$.

The space of horizontal frame-periodic  Euclidean Darboux curves inherits the symplectic structure given by Definition \ref{hypsymp}. This symplectic structure is isomorphic to the structure used by J. Millson and B.A. Zombro in ([17])  as can be easily seen from Lemma \ref{vectprod}. More precisely, in the  Millson-Zombro paper the Euclidean space $E^3$ is identified with $so_3(R)$ which is isomorphic to $su_2$, and their symplectic form is identical to the one given by equation (\ref{dsymp}).

It can be shown by arguments identical to the ones already presented in this paper that the Hamiltonian flow induced by the function $f(x(s))=\half\int_0^l\kappa^2(s)\,ds$ leads to Heisenberg's magnetic equation, and that the passage to the non-linear Schroedinger's equation is the same as the one presented  for the non-Euclidean cases.

The present formalism clarifies   
Hasimoto's first observation that $\psi=\kappa\exp{(i\int\tau\,dx)}$ of a curve $\gamma(s,t)$ that satisfies the filament
equation 
\eq\label{fm}
\frac{\partial \gamma}{\partial t}(s,t)=\kappa (s,t)B(s,t)\ee
is a solution of the non-linear Schroedinger equation. In this notation, it is understood that $t\rightarrow \gamma (s,t)$ denotes a  family of curves in $R^3$ parametrized by $t$ and  that $B(s,t)$ denotes the binormal vector along the curve $s\rightarrow \gamma (s,t)$. 

When the solution curves of the filament equation are restricted to curves  parametrized by arc-length,i.e., to curves $\gamma (s)$ such that $||\frac{d\gamma}{ds}(s)||=1$ then
  \[T(t,s)=
{\partial\gamma\over\partial s}(t,s)\mbox{ and }\frac{dT}{ds}(s)=\kappa (s,t)N(s,t)=\frac{\partial ^2 \gamma}{\partial s^2}(s,t).\]
Moreover, $B(s,t)=T(s,t)\times N(s,t)$. It then follows that in the space of arc-length parametrized curves the filament equation   can be written as
\eq\label{fm1} {\partial\gamma\over\partial t}=
{\partial\gamma\over\partial s}\times{\partial^2\gamma\over
\partial s^2}.\ee
For each solution curve $\gamma(s,t)$ of (\ref{fm1}) the tangent vector $T(s,t)$ satisfies

\eq\label{eucheisen}{\partial T\over\partial t}=T\times{\partial^2T\over
\partial s^2}\ee
as can be easily verified by differentiating with respect to $s$.

 Any solution $T(s,t)$ of the preceeding equation may be interpreted as the coordinate vector of
$\Lambda (s,t)$ relative to an orthonormal basis in the Cartan space $\fp$. Then, $\Lambda (s,t)$ evolves according to\[\frac{\partial\Lambda}{\partial t}(s,t)=i[\frac{\partial^2\Lambda}{\partial s^2}(s,t),\Lambda (s,t)]\] if $\fp=\fh$, and evolves according to \[\frac{\partial\Lambda}{\partial t}(s,t)=[\frac{\partial^2\Lambda}{\partial s^2}(s,t),\Lambda (s,t)]\] if $\fp=\fk$.

The function $\psi (s,t)=u(s,t)\exp{(i\int_0^s u_1(x,t)\,dx)}$ associated with the frame $R(s,t)$ that defines $\Lambda (s,t)$ a solution of the non-linear Schroedinger's equation independently of the choice of the symmetric space (Theorem \ref{h-s}).
 When  the frame $R(s,t)$ is a 
Serret-Frenet frame then:

$u_1=\tau$, $u_2=-\kappa$, $u_3=0$ in the hyperbolic case,

$u_1=\tau+\half$, $u_2=0$, $u_3=\kappa$ in the spherical case, and

$u_1=\tau$, $u_2=0$, $u_3=\kappa$ in the Euclidean case.

In all cases,  \[\psi (s,t)=u(s,t)\exp{(i\int_0^s u_1(x,t)\,dx)}=(\exp{i\theta})(\kappa (s,t)\exp{i\int_0^s\tau (x,t)\,dx)})\] for some angle $\theta$, 
 and therefore, Hasimoto's function $\kappa (s,t)\exp{i\int_0^s\tau (x,t)\,dx)}$ is a solution of the non-linear Schroedinger's equation since the latter is invariant under circular rotations. 
The geometry of the underlying space becomes visible only when the integration of the Hamiltonian equations is carried out  on the full tangent bundle of the Lie group and not just on the part of the equations that resides in the Lie algebra $\fg$.
 \section{ Elastic curves and solitons}
    
    For mechanical systems the Hamiltonian function represents
   the total energy of the system 
and its critical
points  correspond to the equilibrium configurations.  In an  
infinite-dimensional setting the  behaviour of a Hamiltonian system at a critical 
point  of a Hamiltonian seems not to lend itself to such simple characterizations.
 
For the Hamiltonian function $f={1\over 2}\int_0^L k^2\,ds$  the critical points are the  elastic curves.  The solutions of the associated Hamiltonian system that originate at an elastic curve, instead of being stationary, form  travelling waves known as solitons. Soliton solutions of either Heisenberg's magnetic equation or the non-linear Schroedinger's equation are waves that travel at constant speeds with an elastic curve at their wave fronts.
To explain these  statements in some detail it will be necessary to make  a small detour into the geometry of elastic curves.

\subsection{Elastic curves and their Hamiltonian systems}

To maintain continuity with the material already presented and yet to keep the detour at a minimum, the discussion will be confined to the semidirect product $S_K(E^3)$ and $SL_2(C)$. The spherical case $K\times K$ requires adjustments in notation but is otherwise similar to the other two cases (as demonstrated in [11],[12],[13]). 

For notational simplicity $G_\epsilon$  will denote $S_K(E^3)$ for $\epsilon =0$, and $SL_2(C)$ for $\epsilon =-1$. The Lie algebra of $G_\epsilon$ will be denoted by $\fg_\epsilon$. As sets $\fg_0=\fg_{-1}$, but as algebras they are different. Their Lie brackets conform to the following table

\centerline{\begin{tabular}{|c|c|c|c||c|c|c||}
\cline{2-7}
\omit [ , ] & \multicolumn{1}{|c|}{$A_1$} 
	& $A_2$ & $A_3$ & $B_1$ & $B_2$ & $B_3$ \\ \hline
%\omit\hfill[ , ] \hfill\vline& $A_1$ & $A_2$ & $A_3$ 
	%& $B_1$ & $B_2$ & $B_3$ \\ \hline
$A_1$ & $0$ & $-A_3$ & $A_2$ & $0$ & $-B_3$ & $B_2$ \\ \hline
$A_2$ & $A_3$ & $0$ & $-A_1$ & $B_3$ & $0$ & $-B_1$ \\ \hline
$A_3$ & $-A_2$ & $A_1$ & $0$ & $-B_2$ & $B_1$ & $0$ \\ \hline \hline
$B_1$ & $0$ & $-B_3$ & $B_2$ & $0$ & $-\epsilon A_3$ & $\epsilon A_2$ \\ \hline
$B_2$ & $B_3$ & $0$ & $-B_1$ & $\epsilon A_3$ & $0$ & $-\epsilon A_1$ \\ \hline
$B_3$ & $-B_2$ & $B_1$ & $0$ & $-\epsilon A_2$ & $\epsilon A_1$ & $0$ \\ \hline
\end{tabular}}

\centerline{Table 2}

{\definition \label{elast} The problem of finding the minimum of the integral \[\half\int_0^L(u_2^2(s)+u_3^2(s))\,ds\] over all curves $g(s)$ in $G_\epsilon$ that are the solutions of \eq\label{reddarb}\frac{dg}{ds}(s)=g(s)(B_1+u_2(s)A_2+u_3(s)A_3)\ee  and satisfy fixed boundary conditions at $s=0$ and $s=L$ shall be called the elastic problem on $G_\epsilon$.}

{\definition\label{ext} 
 The projections $x(s)=\pi_\epsilon(g(s))$ of the "extremal curves" $g(s)$ on the underlying   space $G_\epsilon/K$ are called elastic curves.}

It is known that  the elastic problem has a solution for any pair of boundary points provided that $L$ is sufficiently large. Since the elastic problem is left-invariant, the initial point
can always be taken at the identity. It is evident from the first part of the paper that \[\half\int_0^L(u_2^2(s)+u_3^2(s))\,ds=\half\int_0^L\kappa^2(s)\,ds\] where $\kappa (s)$ is the curvature of the projected curve $\pi_\epsilon (g(s))$. The set of curves (\ref{reddarb}) may be considered as a "reduced" Darboux space for the function $f(g)=\half\int_0^L\kappa^2(s)\,ds$ for the following reasons:

For any Darboux curve $g(s)$ that is a solution of \[\frac{dg}{ds}(s)=g(s)(B_1+u_1(s)A_1+u_2(s)A_2+u_3(s)A_3)\]
$g_0(s)=g(s)\exp{(-A_1\int_0^su_1(x)\,dx)}$ projects onto the same  base point as $g(s)$ and satisfies (\ref{reddarb}). Consequently,\[f(g(s))=f(g_0(s)).\]

It should be noted that when the solution $R(s)$ of the equation \[\frac{dR}{ds}(s)=R(s)(u_2(s)A_2+u_3(s)A_3)\] is periodic, then $h(s)=g(s)R^*(s)$ is an anchored frame-periodic curve, and therefore, $u(s)=u_2(s)+iu_3(s)$ is a solution of the non-linear Schroedinger's equation  (Theorem \ref{h-s}).
 
  The  solutions of the elastic problem seen from control theoretic perspectives in accordance with the Maximum Principle are confined to the projections of integral curves, called extremal curves, of a Hamiltonian system on the cotangent bundle $T^* {G_\epsilon}$ of $G_\epsilon$ ([11]). 

To take advantage of  the left-invariant symmetries, the cotangent bundle  $T^* G_\epsilon$ will be represented via the left translations as $G_\epsilon \times \fg_\epsilon^*$, where
$\fg_\epsilon^*$ denotes the dual of $\fg_\epsilon$. Linear functions in $\fg_\epsilon^*$ will be represented by the coordinate functions $h_1,h_2,h_3,H_1,H_2,H_3$ relative
to the dual basis $B_1^*,B_2^*,B_3^*,A_1^*,A_2^*,A_3^*$
 defined by the Pauli matrices .

An easy application of the Maximum Principle shows that the regular extremal curves of the 
elastic problem 
are the projections of the integral curves of the Hamiltonian vector field $\vec{H}$ defined by the Hamiltonian  function \eq\label{ham} H=\half(H_2^2+H_3^3)+h_1.\ee The extremal functions $u_2$ and $u_3$ are of the form \eq\label{extctr}u_2=H_2,\,u_3=H_3.\ee
 
{\remark The abnormal extremal curves shall be ignored. 
It is known ([10]) that  the abnormal extremal curves exist, and may or may not project onto an optimal  solutions of the elastic problem. However, the optimal solutions that are the projections of abnormal extremal curves are also the projections of regular extremal curves, and therefore could be ignored, at least as far as optimality is concerned. }

  The most direct way to get the equations of $\vec{H}$ is via the Poisson brackets
involving the variables
$h_1,h_2,h_3,H_1,H_2,H_3$. The Poisson brackets of these variables are isomorphic to the Lie brackets in Table 2, and are reproduced for the convenience of the reader in Table 3 below.
\\[.2in]
\centerline{\begin{tabular}{|c|c|c|c||c|c|c||}
\cline{2-7}
\omit \{ , \} & \multicolumn{1}{|c|}{$H_1$} 
	& $H_2$ & $H_3$ & $h_1$ & $h_2$ & $h_3$ \\ \hline
%\omit\hfill\{ , \} \hfill\vline& $H_1$ & $H_2$ & $H_3$ 
	%& $h_1$ & $h_2$ & $h_3$ \\ \hline
$H_1$ & $0$ & $-H_3$ & $H_2$ & $0$ & $-h_3$ & $h_2$ \\ \hline
$H_2$ & $H_3$ & $0$ & $-H_1$ & $h_3$ & $0$ & $-h_1$ \\ \hline
$H_3$ & $-H_2$ & $H_1$ & $0$ & $-h_2$ & $h_1$ & $0$ \\ \hline \hline
$h_1$ & $0$ & $-h_3$ & $h_2$ & $0$ & $-\epsilon H_3$ & $\epsilon H_2$ \\ \hline
$h_2$ & $h_3$ & $0$ & $-h_1$ & $\epsilon H_3$ & $0$ & $-\epsilon H_1$ \\ \hline
$h_3$ & $-h_2$ & $h_1$ & $0$ & $-\epsilon H_2$ & $\epsilon H_1$ & $0$ \\ \hline
\end{tabular}}
\\[.1in]
\centerline{Table 3}
 Therefore, the Hamiltonian equations are given by:
\begin{equation}\label{hamsys} \begin{array}{rcl}
	{dH_1\over ds}& = &\{H_1,H\}=H_2\{H_1,H_2\}+
		H_3\{H_1,H_3\}+\{H_1,h_1\}=0\\
	{dH_2\over ds}& = &\{H_2,H\}=H_3\{H_2,H_3\}+
		\{H_2,h_1\}=-H_3H_1+h_3\\
	{dH_3\over ds}& = &\{H_3,H\}=H_2\{H_3,H_2\}+
		H_3\{H_3,H_3\}+\{H_3,h_1\}=H_2H_1-h_2\\
	{dh_1\over ds}& = &\{h_1,H\}=H_2\{h_1,H_2\}+
		H_3\{h_1,H_3\}+\{h_1,h_1\}=H_3h_2-H_2h_3\\
	{dh_2\over ds}& = &\{h_2,H\}=H_2\{h_2,H_2\}+
		H_3\{h_2,H_3\}+\{h_2,h_1\}=-H_3h_1+\epsilon H_3\\
	{dh_3\over ds}& = &H_2\{h_3,H_2\}+
		H_3\{h_3,H_3\}+\{h_3,h_1\}=\{h_3,H\}=H_2h_1-\epsilon H_2.
\end{array}\end{equation}

It follows  by an easy calculation that
\[I_1=h_1^2+h_2^2+h_3^2+\epsilon(H_1^2+H_2^2+H_3^2)\] and \[I_2=h_1H_1+h_2H_2+h_3H_3\]are  the constants of motion  for (\ref{hamsys}). Together with $H_1$ and $H$, $I_1$ and $I_2$ account for 
four  independent constants of motion, and therefore, (\ref{hamsys}) is 
completely integrable.

  Any solution $h_1(s),h_2(s),h_3(s),H_1(s),H_2(s),H_3(s)$ of (\ref{hamsys}) defines  complex functions $u(s) =H_2(s)+iH_3(s)$ and $w(s)=h_2(s)+ih_3(s)$. Then, $u(s)$ generates a soliton solution 
if there exists a real number $\xi$ such that $\psi (s,t)=u(s+\xi t)$
is a solution of the non-linear Schroedinger's equation.

It follows from equations (\ref{hamsys}) that
\begin{equation}
\frac{du}{ds}(s)=iH_1u(s)-iw(s) \mbox{, and } \frac{dw}{ds}=i(h_1-\epsilon)u(s).
\end{equation}
Therefore,\[
\frac{\partial\psi}{\partial t}=i\xi (H_1\psi -w)\mbox{, and }\frac{\partial^2\psi}{\partial s^2}
=-H_1^2\psi +H_1w+(h_1-\epsilon)\psi.
\]
Since $H={1\over 2}|\psi |^2 +h_1$,
\begin{eqnarray*}
 -i\frac{\partial\psi}{\partial t}-(\frac{\partial^2\psi}{\partial s^2}+{1\over 2}|\psi |^2\psi)&=&
\xi (H_1\psi -w)-(-H_1^2\psi+H_1w+(h_1-\epsilon)\psi +\psi(H-h_1)) \\&=&
-(\xi +H_1)w+(\xi H_1+H_1^2+\epsilon-H)\psi.
\end{eqnarray*}
The  preceding quantity is a solution of the non-linear Schroedinger's equation whenever\eq\label{sol}
\xi =-H_1\mbox{ and }H=\epsilon.
\ee
{\remark The reader may verify that (\ref{sol}) is a necessary condition  for the existence of solitons for all  three cases. For the spherical case $\epsilon=1$ and $H_1$ and $H$ are the same as in the cases discussed above.}

 To show that periodic solutions $u(s)$ exist on the 
energy level $H=\epsilon$ requires explicit formula for $u(s)$ in terms of the remaining constants of motion $I_1$ and $I_2$.

To begin with, note that 
\[
(H_2h_3-H_3h_2)^2+(H_2h_2+H_3h_3)^2=(H_2^2+H_3^2)(h_2^2+h_3^2).
\]
 Then,
\begin{equation}\label{ellipt}\begin{array}{rcl}
(\frac{d}{ds}h_1)^2 & = & (H_2h_3-H_3h_2)^2 = (H_2^2+H_3^2)(h_2^2+h_3^2)-(H_2h_2+H_3h_3)^2\\
     & =& (H_2^2+H_3^2)(I_1-\epsilon(H_1^2+H_2^2+H_3^2)-h_1^2)-(I_2-h_1H_1)^2\\
     & =& 2(H-h_1)(I_1-\epsilon H_1^2-2\epsilon (H-h_1)-h_1^2)-(I_2-h_1H_1)^2\\
     & =& 2h_1^3+c_1 h_1^2+c_2h_1 +c_3  
\end{array}\end{equation}
where  $c_1,c_2,c_3$ are the constants of motion given by the following expressions
\[c_1=-(H_1^2-2H-4\epsilon), c_2=(2I_2H_1-2\epsilon H_1^2+4\epsilon H-2I_1), c_3= 2H(I_1-\epsilon H_1^2-2\epsilon H)-I_2^2.\]

Therefore, $h_1(s)$ is expressed in terms of elliptic functions, and since 
$k^2=H_2^2+H_3^2 =2(H-h_1)$
 the same can be said  for the curvature of the  projected elastic curve. The remaining 
variables  $u=H_2+iH_3$ and $w=h_2+ih_3$
  can be integrated in terms of  two angles $\theta$ and 
$\phi$ 
  defined as follows:

 \bean 
I_1&=&h_1^2+|w|^2+\epsilon (H_1^2+|u|^2)=h_1^2+|w|^2+\epsilon(H_1^2+2(H-h_1))\\                              &=&  (h_1-\epsilon)^2+|w|^2+\epsilon H_1^2+2\epsilon H-\epsilon ^2\eean
and therefore,
\eq\label{sphere}(h_1-\epsilon)^2+|w|^2=J^2
\ee
where $J^2$ denotes $I_1-\epsilon H_1^2-2\epsilon H+\epsilon^2$. Since $J$ is constant along each extremal trajectory, equation (\ref{sphere}) defines a sphere along each extremal curve. The angles
  $\theta$ and $\phi$  are defined on that sphere by
\begin{equation}\label{polar1}
 (h_1(s)-\epsilon)=J\cos \theta (s)\mbox{ and } w(s)=J\sin \theta (s)e^{i\phi (s)}.
\end{equation}
 It follows that 
\eq\label{polar2}
\frac{dh_1}{ds}=-J\sin\theta\frac{d\theta}{ds} \mbox{, and  }\frac{dw}{ds}=
w(\frac{cos\theta}{\sin \theta}\frac{d\theta}{ds}+i\frac{d\phi}{ds}).
\ee
 Furthermore,
\begin{eqnarray*}
\frac{u}{w}=\frac{u\bar{w}}{|w|^2}&=& \frac{H_2h_2+H_3h_3+i(H_3h_2-H_2h_3)}{J^2-(h_1+\epsilon)^2}\\
                                  &=& \frac{I_2-h_1H_1+i\frac{dh_1}{ds}}{J^2\sin ^2\theta}\\
              &=&\frac{I_2-h_1H_1}{J^2\sin ^2\theta}-{i\over J\sin \theta}\frac{d\theta }{ds}.
\end{eqnarray*}
Therefore,
\begin{eqnarray*}
w(\frac{\cos\theta}{\sin\theta}\frac{d\theta }{ds}+i\frac{d\phi}{ds})&=&\frac{dw}{ds}\\
&=&i(h_1-\epsilon )u  = (iJ\cos\theta \frac{(I_2-h_1H_1)}{J^2\sin ^2\theta}+
\frac{\cos\theta}{\sin\theta}\frac {d\theta}{ds})w
\end{eqnarray*}
hence,
\begin{equation}
\frac{d\phi}{ds}=J\cos\theta \frac{(I_2-h_1H_1)}{J^2\sin ^2\theta}=\frac{J\cos\theta(I_2-\epsilon H_1-H_1J\cos\theta)}{J^2\sin ^2\theta}.
\end{equation}
Equation (\ref{ellipt}) can be also written as
\[(\frac{dh_1}{ds})^2=2(H-h_1)|w|^2-(I_2-H_1h_1)^2.\]
 The substitutions from (\ref{polar1}) and (\ref{polar2}) in the  preceeding equation define  $\theta$ as the solution of the following differential equation
\begin{equation}
{(\frac{d\theta}{ds})}^2=2(H-\epsilon-J\cos \theta )-\frac{(I_2-H_1(\epsilon+J\cos\theta))^2}{J^2\sin^2\theta}.
\end{equation}
It follows from  (\ref{polar1}) and the equation that follows that both $u(s)$ and $w(s)$ are
determined by  the angles $\phi (s)$ and $\theta (s)$.
We now return to the question of periodicity of $u(s)$. Evidently, both $u$ and $w$ are periodic  whenever $\phi (0)=\phi (L)$ and $\theta (0)=\theta (L)$. 
 Soliton solutions propagate with speed $\psi =H_1$ on the energy level $H=\epsilon$. On this energy level   $\phi (0)=\phi (L)$ and $\theta (0)=\theta (L)$ if and only if
\begin{equation}
\int_0^L\frac{J\cos\theta(I_2+H_1-H_1J\cos\theta)}{J^2\sin ^2\theta}\,ds=0
\end{equation}
 where $\theta$ denotes a closed solution of the equation
\begin{equation}
{(\frac{d\theta}{ds})}^2=-2J\cos \theta -\frac{(I_2-H_1(\epsilon+J\cos\theta))^2}{\sin^2\theta}.
\end{equation}

It is known that there are are infinitely many closed solutions
for suitable constants $I_1,I_2,H_1$ ( for instance, ([9])) however, such calculations will not be presented here.  

\section{Complete Integrability}

There are further connections between elasic curves and solutions of the non-linear 
Schroedinger equation
that were first noticed by J. Langer and R. Perline, namely that some of the integrals  of motion for the elastic curves correspond to the integrals of motion for the non-linear Schroedinger's equation ([14]).  We will illustrate this phenomenon by showing that the
function $f(\gamma )=\int_0^L k^2\tau\,ds$ is an integral of motion for Heisenberg's magnetic 
equation,
 while the quantity $k^2(s)\tau (s)$ is a constant of motion for the elastic 
problem.

 We shall first show that $k^2\tau$ is a constant of motion for the elastic curves.
 The Serret-Frenet frame and the frame in Definition \ref{elast}  
rotate around each other in the plane perpendicular to the tangent vector. If we denote by 
$\beta$
the angle through which the frame in Definition \ref{elast} rotates relative to the Serret-Frenet frame then, according to the formula in ([12], p461)
  \[\frac{d\beta}{ds}(s)=\tau\mbox{ and }\tan\beta= -{H_2\over H_3}.\]  

It follows that
\begin{eqnarray*}
\sec ^2\beta\frac{d\beta}{ds}&=&\frac{(H_2\frac{dH_3}{ds}-H_3\frac{dH_2}{ds})}{H_3^2}\\
                             &=&\frac{H_2(H_2H_1-h_2)-H_3(-H_3H_1+h_3)}{H_3^2}\\
                             &=&\frac{H_1(H_2^2+H_3^2)-H_2h_2-H_3h_3}{H_3^2}\\
                             &=&\frac{H_1(H_2^2+H_3^2)-(I_2-H_1h_1)}{H_3^2}\\
                             &=&\frac{2H_1H-I_2}{H_3^2}.
\end{eqnarray*}
  Since \[\sec^2\beta=(H_2^2+H_3^2)/H_3^2=\kappa^2/H_3^2\] it follows that
\[\kappa ^2\tau =2H_1H-I_2=constant.\]

{\theorem  (a) The Hamiltonian flow of $f_1=\int_0^l\kappa^2 (s)\tau (s)\,ds$ is given by
\eq\label{ksqtau}
\frac{\partial \Lambda}{\partial t}=2(\dddot\Lambda -\leangle \dddot{\Lambda},\Lambda\riangle\Lambda )-3\leangle \Lambda,
\ddot{\Lambda}\riangle\dot{\Lambda}).
\ee

(b) $f_0=\half\int_0^L\kappa^2(s)\,ds$ and $f_1$ Poisson commute. }

\prf 

The first part of the proof consists in showing that \[\kappa^2\tau =-i\leangle [\Lambda ,\frac{d\Lambda}{ds}],\frac{d^2\Lambda}{ds^2}\riangle.\]

Suppose now that $T(s)=\Lambda (s)$ denotes the Hermitian matrix that corresponds to the tangent 
vector of a horizontal Darboux curve that projects onto a curve
$\gamma \in \HI^3$. Then $N(s)$ and $B(s)$, the matrices that correspond to the normal and the binormal vectors, are given by \[N={1\over \kappa}
\frac{d\Lambda}{ds}\mbox{  and }B(s)={1\over i}[T(s),N(s)]=-\frac{i}{\kappa}[\Lambda ,\frac{d\Lambda}{ds}].\]
 According to 
the Serret-Frenet
equations $\frac{dN}{ds} = -k \Lambda + \tau B $. Therefore, \[\tau =\leangle \frac{dN}{ds},
B \riangle =
-i\leangle -{1\over k^2}\frac{dk}{ds}\frac{d\Lambda}{ds}+{1\over k}\frac{d^2\Lambda}{ds^2}, 
{1\over k}[\Lambda,
\frac{d\Lambda}{ds}]\riangle .\]
It follows that \[k^2\tau =-i\leangle [\Lambda ,\frac{d\Lambda}{ds}],\frac{d^2\Lambda}{ds^2}\riangle.
\] 
Let $V(s)$ be an arbitrary tangent vector at a  frame-periodic horizontal Darboux curve $g(s)$. Then the directional derivative of $f_1$ at 
$g$ in the 
direction $V$ is given by the following expression:

 \[df_1(V) = -i\frac{\partial}{\partial t}\int_0^L\langle [Z(s,t),\frac{\partial Z}{\partial s}(s,t)],\frac{\partial ^2 Z}{\partial s^2}(s,t)\rangle \,ds|_{t=0}\]
where the dots indicate derivatives with respect to $s$, and where $Z(s,t)$ denotes a field of Hermitian matrices such that \[Z(s,0=\Lambda (s)\mbox{ and }\frac{\partial Z}{\partial t}(s,0)=V(s).\] 

It follows that
\begin{eqnarray*} df_1(V) &=&-i\int_0^L\leangle \dddot{V}, [\Lambda,\dot {\Lambda}]\riangle+
\leangle \ddot{\Lambda} ,[\dot{V},
\dot{\Lambda} ]\riangle+\leangle \ddot{\Lambda},[\Lambda,\ddot{V}]\riangle\,ds\\
&=& -i\int_0^L 2\leangle [\ddot{\Lambda} ,\Lambda ],\ddot{V}\riangle-\leangle [\ddot{\Lambda} ,
\dot{\Lambda}],
\dot{V}\riangle\,ds\\
    &=& i\int_0^L \leangle 2(\frac{d}{ds}([\ddot{\Lambda},\Lambda ])+[\ddot{\Lambda},
\dot{\Lambda}]),\dot{V}\riangle\,ds\\
&=&i\int_0^L\leangle 2[\dddot{\Lambda},\Lambda ])+3[\ddot{\Lambda},
\dot{\Lambda}],\dot{V}\riangle\,ds.
\end{eqnarray*}

Let $V_1(s)$ denote the Hermitian matrix such that\[df_1(V)=\omega_{\Lambda}(V_1,V).\]

Then,
\[i\int_0^L\leangle 2[\dddot{\Lambda},\Lambda ]+3[\ddot{\Lambda},
\dot{\Lambda}],\dot{V}\riangle\,ds={1\over i}\int_0^L\langle [\Lambda ,\dot{V_1}],\dot{V}\rangle\,ds\] which implies \[
\int_0^L\langle [\Lambda ,\dot{V_1}]+2[\dddot{\Lambda},\Lambda ]+3[\ddot{\Lambda},\dot{\Lambda}],\dot{V}\rangle\,ds=0.
\]
 When $\dot{V}=[\Lambda,C]$ the above becomes\[\int_0^L\langle [[\Lambda ,\dot{V_1}],\Lambda ]+2[[\dddot{\Lambda},\Lambda ],\Lambda ]+3[[\ddot{\Lambda},\dot{\Lambda}],\Lambda ],C(s)\rangle\,ds=0.
\]
Since $C(s)$ can be an arbitrary curve with $C(0)=0$  the preceeding integral equality reduces to 
\[
[[\Lambda ,\dot{V_1}],\Lambda ]+2[[\dddot{\Lambda},\Lambda ],\Lambda ])+3[[\ddot{\Lambda},\dot{\Lambda}],\Lambda ]=0.
\]
 The Lie bracket relations in Lemma \ref{triplebracket} imply that
\[
\dot{V_1}+2(\langle \dddot{\Lambda} ,\Lambda\rangle\Lambda-\dddot{\Lambda})+3\langle \ddot{\Lambda} ,\Lambda\rangle \dot{\Lambda}=0.
\]
 Now it follows by the arguments  used earlier in the paper that the Hamiltonian flow $X_{f_1}$ satisfies\[\frac{\partial \Lambda}{\partial t}=2(\dddot\Lambda -\leangle \dddot{\Lambda},\Lambda\riangle\Lambda )-3\leangle \Lambda,
\ddot{\Lambda}\riangle\dot{\Lambda}).\]
 Thus part (a) is proved. 

To prove  part (b) it is required to show that
  the Poisson bracket  of $f_0$ and $f_1$, given by the formula
 \[\{f_0,f_1\}(\Lambda)=\omega_\Lambda (V_0(\Lambda), V_1(\Lambda))={1\over i}\int_0^L\langle \Lambda (s),[\dot{V_0(s)},\dot{V_1}(s)]\,ds\]
with $\dot{V}_0(\Lambda)=i[\ddot\Lambda ,\Lambda]$ and $\dot{V_1}=-(2(\langle \dddot{\Lambda} ,\Lambda\rangle\Lambda-\dddot{\Lambda})+3\langle \ddot{\Lambda} ,\Lambda\rangle \dot{\Lambda})$, is equal to $0$.

An easy calculation shows that
\[
[\dot{V_0},\dot{V_1}]=i(2(\leangle\dddot{\Lambda} ,\ddot\Lambda\riangle)-\leangle \Lambda ,\ddot {\Lambda }\riangle
\leangle \dddot {\Lambda} ,\Lambda \riangle) -3\leangle \ddot {\Lambda }
,\Lambda \riangle
\leangle \ddot {\Lambda} ,\dot{\Lambda }\riangle ) 
\Lambda .
\]
 
Hence,\[\{f_0,f_1\}=\int_0^L (2\leangle\dddot\Lambda ,\ddot\Lambda\riangle -
2\leangle \Lambda ,\dddot\Lambda\riangle\leangle \Lambda ,\ddot\Lambda\riangle
 -3\leangle\Lambda ,\ddot\Lambda\riangle\leangle\dot\Lambda ,\ddot\Lambda\riangle)\,ds.
\]
 
The integral of the first term is zero because $2\leangle \dddot\Lambda ,\ddot\Lambda\riangle=
\frac{d}{ds}\leangle \ddot\Lambda ,\ddot\Lambda\riangle$.

 Since 
$2\leangle\Lambda ,\ddot\Lambda\riangle\leangle \Lambda ,\dddot\Lambda\riangle =
\frac{d}{ds}\leangle \Lambda ,\ddot\Lambda\riangle^2 -2\leangle \Lambda ,\ddot\Lambda
\riangle\leangle 
\dot\Lambda ,\ddot\Lambda\riangle$, the remaining integrand reduces to one term $-\leangle\Lambda ,
\ddot\Lambda\riangle
\leangle \dot\Lambda ,\ddot\Lambda\riangle$. But then 
${1\over 4}\frac{d}{ds}\leangle\dot\Lambda ,\dot\Lambda\riangle^2 =\leangle \Lambda ,
\ddot\Lambda\riangle\leangle\dot\Lambda ,\ddot\Lambda\riangle $ because
$\leangle\dot\Lambda ,\dot\Lambda\riangle=-\leangle\Lambda ,\ddot\Lambda\riangle$, and 
part (b) is proved.
\qed
{\theorem Suppose that $\Lambda (s,t)=R(s,t)B_1R^*(s,t)$ evolves according to  the equation (\ref{ksqtau}) where $R(s,t)$ is the solution of \[\frac{\partial R}{\partial s}(s,t)=R(s,t)\bmt 0&u(s,t)\\-\bar{u}(s,t)&0\emt ,\,R(0,t)=I.\] Then, $u(s,t)$ is a solution of \eq\label{kdv}\frac{\partial u}{\partial t}-3|u|^2\frac{\partial u}{\partial s}-2\frac{\partial ^3 u}{\partial s^3}=0.\ee}
 
 This theorem is proved by a calculation similar to the one used in the proof of Theorem \ref{h-s}, the details of which will not be reproduced here.

 Equation (\ref{kdv}) is similar to the modified Korteweg-de Vries equation ( R.Abraham and J. Marsden ([1])) 
   \eq\label{modkdv}v_t-6v^2 v_x+v_{xxx}=0\ee with some notable differences.  Equation (\ref{kdv}) is a complex equation while the modified Korteweg-de Vries equation is a real equation. Because of the difference in sign in front of the third derivative it is not apparent that equation (\ref{kdv}) is the complexification
(modulo some homothetical transformation) of the equation (\ref{modkdv}). It remains an intriguing question if there are any connections between the Korteweg-de Vries equation and the elastic curves.

Functions $f_0$ and $f_1$ also appear in a paper on integrability of the non-linear Schroedinger's equation by C. Shabat and V. Zacharov ([18]), as first noticed by J. Langer and R. Perline ([14]), but in a completely different context. The first two integrals of motion in the paper of Shabat and Zacharov are up to the constant factors given by the following integrals:
\[C_1=\int_{-\infty}^{\infty}|u(s,t)|^2\,ds,\,C_2=\int_{-\infty}^{\infty}(u(s,t))\dot{\bar{u}}(s,t))-\bar{u}(s,t)\dot{u}(s,t))\,ds\]
where they are interpreted as the number of particles and the momentum. To see that $C_1$ and $C_2$ are in exact correspondence with functions $f_0$ and $f_1$  assume that the Darboux curves are expressed by reduced frames $R(s)$ as in Definition \ref{elast} i.e., as the solutions of\[\frac{dR}{ds}(s)=R(s)U(s)\mbox { with } U(s)=u_2(s)A_2+u_3(s)A_3\]
Then,
\[f_0=\half\int_0^L||\dot{\Lambda}(s)||^2\,ds=\half\int_0^L||[B_1,U(s)]||^2,ds=\half\int_0^L|u(s)|^2\,ds.\]
Hence, $ C_1$ corresponds to $\int_0^L\kappa^2(s)\,ds$. Furhermore, $f_1=\int_0^L k^2\tau\,ds$ can be written as \[f_1={1\over 2i}\int_0^L(u(s)\dot{\bar{u}}(s)-\bar{u}(s)\dot{u}(s))\,ds\]because \[
i\leangle \Lambda ,[\dot\Lambda ,\ddot\Lambda]\riangle=\leangle [[B_1,U],[B_1,\dot U]],B_1 \riangle=
\it{Im}\bar{ u}\dot u,\] and  therefore  $f_1$ corresponds to $C_2$.
 
  In the language of mathematical physics the vector $\int_0^L\Lambda (s)\,ds$ is called the total spin ([5]). In this paper 
it appears
as the moment map discussed in the previous section. It is a conserved quantity since the 
Hamiltonian is 
invariant under the action of $SU_2$. This fact can be verified directly as follows:
\[
\frac{\partial}{\partial t}\int_0^L\Lambda(t,s)\,ds =\int_0^L\frac{\partial\Lambda }
{\partial t}(t,s)\,ds=
i\int_0^L[\frac{\partial^2\Lambda}{\partial s^2},\Lambda]\,ds=i\int_0^L\frac{\partial}
{\partial s}[\Lambda ,
\dot\Lambda]\,ds=0.
\]
The third integral of motion  $C_3$ in ([18]), called the energy, is given by \[C_3=\int_{-\infty}^{\infty}(|\frac{\partial u}{\partial s}(s,t)|^2 -\frac{1}{4}|u(s,t)|^4)\,ds.\]
  It corresponds to the function \[f_2=\int_0^L(||\ddot{\Lambda}(s)||^2-\frac{5}{4}||\dot{\Lambda}(s)||^4 )\,ds=\int_0^L (\frac{\partial \kappa}{\partial s}(s)^2+\kappa^2(s)\tau^2(s)-\frac{1}{4}\kappa^4(s) )\,ds.\]
 Functions $f_0,f_1,f_2$ are in involution, i.e., they Poisson commute pairwise. There is a hierarchy of functions that contains $f_0,f_1,f_2$  such that any two functions in the hierarchy Poisson commute. For instance, D.Krepski has shown that $f_3=\int_0^L \tau (s)\,ds$ is in this hierarchy, and he has also shown that the flow of the corresponding Hamiltonian vector field generates the curve shortening equation ([4])\[\frac{\partial\Lambda }{\partial t}(s,t)=\frac{\partial \Lambda}{\partial s}(s,t)=\kappa (s,t)N(s,t).\]

A detailed investigation of  this hierarchy of Poisson commuting functions and its relation to the hierarchies obtained by Langer-Perline and Shabat-Zacharov will be deferred to a separate study.

\begin{center}
\textbf{References}
\end{center}

[1]. Abraham R., and Marsden J., \textit{Foundations of Mechanics},
 Benjamin- Cummings( 1978),
 Reading,Mass

[2]. Arnold V.I., and Khesin B.A., \textit{Topological Methods in Hydrodynamics},  
App. Math. Sci.(125), Springer-Verlag (1998), New York

[3]. Brylinski J.P., \textit{Loop Spaces, Characteristic Classes and Geometric Quantization},
 Progress in Math.(108), Birkhäuser(1993), Boston

[4]. Epstein C.L., and Weinstein M.I., A stable manifold theorem for the curve shortening equation,
 Comm. Pure and App. Math., Vol XL (1987),p 119-139

[5]. Faddeev L., and Takhtajan L., \textit{ Hamiltonian Methods in the Theory of Solitons},
 Springer-Verlag (1980), Berlin

[6]. Hamilton R.S., The inverse function theorem of Nash and Moser, Bull. Amer. Math. Soc.
(7), (1972), p 65-221

[7]. Hasimoto H., Motion of a vortex filament and its relation to elastica, J. Phys. Soc. Japan,
(31), (1971), p 293-298

[8]. Hasimoto H., A soliton on a vortex filament, J. Fluid Mech. (51),(1972), p 477-485

[9]. Ivey T., and Singer D.A., Knot types, homotopies and stability of closed elastic curves, 
Proc. London Math.Soc (3) 79, (1999),pp 429-450

[10]. Jurdjevic V., Integrable Hamiltonian systems on Complex Lie groups, a preprint

[11]. Jurdjevic V..\textit{Geometric Control Theory}, Cambridge Studies in Advanced Mathematics(51), Cambridge Univ. Press (1997), New York

[12]. Jurdjevic V., and Monroy-Perez F., Hamiltonian systems on Lie groups: Elastic curves,
 Tops and Constrained Geodesic Problems, \textit{Non-Linear Geometric Control Theory and
 its Applications}, World Scientific Publishing Co., (2002), Singapore

[13]. Jurdjevic V., Hamiltonian Systems on Lie Groups: Kowalewski type, Ann. Math.,(150),
(1999), p 1-40

[14]. Langer J., and Perline R., Poisson Geometry of the Filament Equation,
 J. Nonlinear Sci., Vol 1, (1978), p 71-93 

[16]. Magri F. A simple model for the integrable Hamiltonian equation, J. Math.
Phys. (19), (1978), p 1156-1162

[17]. Millson J., and Zombro B.A., A K\"ahler structure on the moduli spaces of isometric
maps of a circle into Euclidean spaces, Invent. Math. Vol 123, (1), (1996), p35-59

[18]. Shabat C., and Zakharov V., Exact theory of two dimensional self-focusing
 and one dimensional self-modulation of waves in non-linear media, Sov. Phys. JETP, (34),(1972), p 62-69

[19]. Sternberg S., \textit{Lectures on Differential Geometry}, Prentice-Hall Inc.(1964),
 Englewood-Cliffs, New Jersey

[20]. Sulem C., and Sulem P-L., \textit{The nonlinear Schroedinger Equation; Self Focusing and Collapse},
 Springer-Verlag (1999), New York

[21]. Wagneur B., The symplectic structure of loops in 3-dimensional hyperbolic spaces, PhD thesis(2000),
The University of Maryland, Maryland

\end{document}